\documentclass[journal ]{new-aiaa}
\usepackage[utf8]{inputenc}
\usepackage{textcomp}
\usepackage{subcaption}
\usepackage{graphicx}
\usepackage{amsmath}
\usepackage[version=4]{mhchem}
\usepackage{siunitx}
\usepackage{longtable,tabularx}
\usepackage{bm}
\usepackage{tikz}
\usetikzlibrary{calc,angles,positioning,intersections,decorations.pathreplacing,decorations.markings}
\DeclareMathOperator{\diag}{diag}

\usepackage{amsthm}
\newtheorem*{remark}{Remark}
\newtheorem{conjecture}{Conjecture}

\title{Constrained Fuel and Time Optimal 6DOF Powered Descent Guidance Using Indirect Optimization}

\author{Nicholas P. Nurre \footnote{Graduate Research Assistant, Department of Aerospace Engineering, 141 Engineering Dr, AIAA Student Member.} and Ehsan Taheri\footnote{Assistant Professor, Department of Aerospace Engineering, 141 Engineering Dr, AIAA Senior Member.}}
\affil{Auburn University, Auburn, AL, 36849}

\begin{document}

\maketitle

\begin{abstract}
Powered descent guidance (PDG) problems subject to six-degrees-of-freedom (6DOF) dynamics allow for enforcement of practical attitude constraints. However, numerical solutions to 6DOF PDG problems are challenging due to fast rotational dynamics coupled with translational dynamics, and the presence of highly nonlinear state/control path inequality constraints. In this work, constrained fuel- and time-optimal 6DOF PDG problems are solved leveraging a regularized indirect method, subject to inequality constraints on the thrust magnitude, thruster gimbal angle, rocket tilt angle, glideslope angle, and angular velocity magnitude. To overcome the challenges associated with solving the resulting multipoint boundary-value problems (MPBVPs), the state-only path inequality constraints (SOPICs) are enforced through an interior penalty function method, which embeds the resulting MPBVPs into a multi-parameter smooth neighboring families of two-point BVPs. Extremal solutions are obtained using an indirect multiple-shooting solution method with numerical continuation. Moreover, an empirical relation is derived for the directly-adjoined Lagrange multipliers associated with SOPICs. The  fuel- and time-optimal trajectories are compared against solutions of DIDO --- a capable pseudospectral-based software for solving practical constrained optimal control problems.   
\end{abstract}

\section{Introduction}
\lettrine{P}{owered} descent guidance (PDG) has been well studied since the 1960s \cite{meditch_problem_1964} and was first implemented when NASA used an explicit guidance law in the Apollo missions \cite{klumpp_apollo_1974}. Nowadays, computational-based guidance methods represent the state-of-the-art \cite{lu_introducing_2017} and have enabled, most notably, SpaceX to autonomously and safely perform extremely complicated rocket landings \cite{blackmore_autonomous_2016} with convex optimization \cite{malyuta_convex_2021,wang_survey_2024}. Unlike three-degree-of-freedom (3DOF) PDG problems, the 6DOF PDG problem allows for the enforcement of non-trivial attitude dynamics constraints \cite{malyuta_advances_2021}. A difficulty, however, lies in the presence of highly nonlinear couplings between translational/rotational dynamics and existence of various control/state constraints, making 6DOF PDG problems quite challenging to solve. However, direct optimization methods \cite{trelat_optimal_2012} have been shown to provide reliable solutions. The recent computational advances in convex optimization allow this class of PDG problems to be solved in real time \cite{scharf_implementation_2017}. For instance, Refs.~\cite{szmuk_successive_2016,szmuk_successive_2018} solve the 6DOF PDG problem with a convex-based approach, which is extended to handle state-triggered constraints in Ref. \cite{szmuk_successive_2020}. Ref. \cite{sagliano_six-degree--freedom_2024} offers another convex-based approach to the problem using augmented
convex–concave decomposition. 

Direct methods are more practical in tackling 6DOF PDG problems, but may require denser grid points to generate accurate solutions due to discretization errors. They also may have difficulty in solving optimal control problems (OCPs) with singular control arcs \cite{mall_solving_2021,andres-martinez_switched_2022,pager_method_2022} and can't offer concrete insights into the control structure of the problem. As such, Refs. \cite{reynolds_optimal_2020,malyuta_advances_2021} point out the optimal solution to the 6DOF PDG problem remains an open problem. Indirect-based methods offer high-resolution solutions that can provide greater accuracy. It is also possible to obtain entirely new classes of solutions to OCPs which were once believed to have only a particular extremal control profile. For instance, Ref.~\cite{mall_using_2023} used the Uniform Trigonometrization Method (UTM) \cite{mall_uniform_2020} to solve a well-known free-flyer problem \cite{sakawa_trajectory_1999}. The authors present a completely new solution that consists of a singular control arc, outperforming the bang-bang solutions reported by Betts \cite{betts_optimal_2003} (see Chap. 6) and also the solution obtained with GPOPS-II \cite{patterson_gpops-ii_2014}. Indirect methods can also be used for generating large databases of optimal and near-optimal solutions \cite{izzo_real-time_2021,zhang_efficient_2024}. For these reasons, indirect-based approaches continue to retain their importance, in particular, to determine the theoretically guaranteed extremal solutions. In fact, explicit guidance laws such as the Apollo guidance is proven to be optimal (under certain assumptions) through the application of the indirect method. A recent example of the theoretical insights gained by indirect methods is Ref. \cite{you_theoretical_2022} where fuel-optimal 3DOF PDG problems subject to a glideslope angle constraint are studied and it is proved that 1) the optimal thrust magnitude has a bang-bang profile (i.e., no singular control arc), 2) the optimal thrust magnitude can have more than three subarcs, and 3) the glideslope constraint can only be active at isolated points. These theoretical insights into the structure of the optimal solution (e.g., the number of control switches and the time of switches) allows for identifying a number of key features in OCPs and leveraging them for developing new algorithms, as is achieved in \cite{kenny_optimal_2024}. Indirect approaches involve deriving the necessary conditions of optimality and numerically solving the resulting boundary-value problem, either through root-finding or collocation. Our studies indicate that solving OCPs with different solution methodologies is illuminating and enriches our understanding of the potential solutions. Our goal, however, is not to compare advantages and disadvantages of \textit{direct} and \textit{indirect} methods, which are already addressed comprehensively in the literature \cite{betts_optimal_2003}, but to advance the state of the art in solving 6DOF PDG problems using enhanced indirect methods.

In the literature, translational-only 3DOF PDG problems are studied in Refs.~\cite{lu_propellant-optimal_2018,lu_propellant-optimal_2023,leparoux_structure_2022,spada_directindirect_2023} using the indirect methods. Ref.~\cite{lu_propellant-optimal_2018} provides a proof for the thrust magnitude switching structure of the 3DOF PDG problem and develops an analytical guidance law based on the structure of the thrust profile. This work is extended in Ref.~\cite{lu_propellant-optimal_2023} by considering the thrust pointing control constraint that leads to approximate open-loop control expressions. The bang-bang thrust control is regularized using the hyperbolic tangent smoothing \cite{taheri_generic_2018}. Ref.~\cite{leparoux_structure_2022} derives exact open-loop control expressions for the 3DOF PDG problems under a thrust pointing control constraint and a glideslope constraint. Ref.~\cite{spada_directindirect_2023} develops an indirect method that can rapidly generate solutions to the 3DOF PDG problems. Refs.~\cite{reynolds_optimal_2020,kovryzhenko_generalized_2024} investigate the 3DOF PDG problem with 2 translational and 1 rotational degrees of freedom, which is useful for considering simple attitude constraints. Ref.~\cite{reynolds_optimal_2020} solves the problem by deriving the exact open-loop control expressions and also provides a proof on the control structure. Ref. \cite{kovryzhenko_generalized_2024} solved the problem through a generalized vectorized trigonometrization method, which is an advanced indirect method for solving constrained OCPs. 

The presence of different state-path inequality constraints is another challenging aspect of 6DOF PDG problems, since they can alter the structure of optimal control profiles. Unlike the direct methods, the indirect method still heavily lacks in practical enforcement methods for these types of constraints. Two forms of the necessary conditions that have been derived are the direct and indirect adjoining approaches \cite{hartl_survey_1995}. Practical implementation of these necessary conditions using traditional indirect methods has been difficult because the sequence of active and inactive constraints along with their time duration must be known \textit{a priori}. They can also introduce corner and/or jump discontinuities in the costates/Hamiltonian. These conditions have to be derived depending on the order of the constraints and result in notoriously difficult-to-solve MPBVPs. However, in recent years, notable progresses are made in detecting switch times \cite{aghaee_switch_2021} and in overcoming the challenges associated with incorporating state-only path inequality constraints (SOPICs) \cite{wang_indirect_2017,antony_path_2018,mall_uniform_2020,mall_solving_2021} and offer promising avenues for solving practical/challenging optimal control problems \cite{epenoy_fuel_2011,perez-palau_fuel_2018}. The three main contributions are 1) formulation and solution to the fuel- and time-optimal 6DOF PDG problems with SOPICs using the indirect method, 2) derivation of the optimal closed-form thrust and thrust steering control expressions under a gimbal-angle constraint, and 3) an approximate relation for the Lagrange multipliers associated with the SOPICs when directly adjoined to the Lagrange cost. To our best knowledge, this is the first application of indirect methods for solving the 6DOF PDG problem. We consider inequality constraints on the thrust magnitude, thrust gimbal angle, rocket tilt angle, angular velocity magnitude, and glideslope angle. Control constraints are enforced by directly adjoining them to the Hamiltonian, allowing closed-form expressions to be derived. The resulting expressions are piecewise continuous and are regularized for numerical implementation. These control expressions are not closed-form under the SOPICs. The SOPICs are enforced with an interior penalty function method \cite{kovryzhenko_generalized_2024}. The Lagrange multipliers associated with the SOPICs, when directly adjoined to the Lagrange cost, are found based on the Hamiltonian invariance principle \cite{taheri_costate_2021}. This relation allows us to verify satisfaction of complementarity slackness \textit{a posteriori}, providing additional checks on the optimality of the resulting solutions. Despite presenting it as a conjecture without rigorous proof, we empirically validate this relation by solving the fuel- and time-optimal 6DOF PDG problems with DIDO \cite{ross_enhancements_2020} -- a capable pseudospectral-based software for solving practical OCPs -- and comparing with our indirect results. We also empirically validate the conjecture by solving the Breakwell problem which has an analytic solution. 

The remainder of the paper is organized as follows. The dynamical model of the 6DOF PDG problem are presented in Section \ref{sec: dynamical model}. Section \ref{sec: optimal control problem} outlines the formulation of fuel- and time-optimal trajectory optimization problems. Section \ref{sec: derivation of necessary conditions} derives the necessary conditions of optimality for both fuel- and time-optimal problems. Section \ref{sec: boundary value problem} describes how the resulting OCPs are solved. Section \ref{sec: results} presents the results. Finally, Section \ref{sec: conclusion} concludes the paper.

\section{Dynamical Model} \label{sec: dynamical model}

The dynamics and kinematics are modeled according to Ref.~\cite{szmuk_successive_2016} with the addition of the aerodynamics terms taken from \cite{sagliano_six-degree--freedom_2024}. Similar to Ref.~\cite{szmuk_successive_2016} and unlike Ref.~\cite{sagliano_six-degree--freedom_2024}, we do not consider a mechanism to control the rocket's roll angle. While our model is 6DOF, it is only controllable in 5 of them. However, our model may be extended to include roll control either through reaction control thrusters introduced in Ref.~\cite{sagliano_six-degree--freedom_2024}, aerodynamic deflection surfaces, or by considering multiple engines that are offset from the rocket's longitudinal axis. 

\begin{figure}
    \centering
    \begin{tikzpicture}

        \begin{scope}[scale=0.75]
            \draw[black,->] (0,0) -- (2,0) node[pos=1.12] {\(\hat{\bm{y}}_\mathcal{I}\)};
            \draw[black,-,dashed] (0,0) -- (-2,0);
            \draw[black,->] (0,0) -- (0,2) node[pos=1.12] {\(\hat{\bm{z}}_\mathcal{I}\)};
            \draw[black,->] (0,0) -- (-1,-1) node[pos=1.12] {\(\hat{\bm{x}}_\mathcal{I}\)};

            \filldraw[black] (0,0) circle (0.05);

            \draw[purple,rotate=70] (0,0) -- (0,1.5);
            \draw[purple,rotate=-70] (0,0) -- (0,1.5);
            \draw[purple] (0,1.5*cos{70}*1.05) ellipse ({1.5*sin{70}} and {1.5*sin{70}*0.1});

            \draw[black,dashed,->] ([shift=(0:1)]0,0) arc (0:20:1)  node[pos=0,black,below] {\(\gamma_\text{min}\)};

        \end{scope}

        \def\xB{-5}
        \def\yB{4}

        \begin{scope}[shift={(\xB,\yB)},scale=4]
            \draw[blue,dashed] (0,0) -- (1,0);
            \draw[blue,dashed] (0,0) -- (0,1);
            \draw[blue,dashed] (0,0) -- (-0.5,-0.5);

            \draw[black,dashed,->] ([shift=(90:0.8)]0,0) arc (90:0:0.8)  node[pos=1.05,black] {\(\theta_\text{max}\)};

            \draw[black,dashed,->] ([shift=(90:0.9)]0,0) arc (90:80:0.9)  node[pos=1.3,black] {\(\theta\)};
                    
            \begin{scope}[rotate around={-10:(0,0)}]
                \draw[blue,->] (0,0) -- (1,0) node[pos=1.12,black] {\(\hat{\bm{y}}_\mathcal{B}\)};
                \draw[blue,->] (0,0) -- (0,0) -- (0,1) node[pos=1.12,black] {\(\hat{\bm{z}}_\mathcal{B}\)};
                \draw[blue,->] (0,0) -- (0,0) -- (-0.5,-0.5) node[pos=1.12,black] {\(\hat{\bm{x}}_\mathcal{B}\)};

                \filldraw[blue] (0,0) circle (0.025);

                \def\armLength{0.5}
                \begin{scope}[shift={(0,-\armLength)}]
                    \draw[green,rotate=20] (0,0) -- (0,1);
                    \draw[green,rotate=-20] (0,0) -- (0,1);
                    \draw[green] (0,0.94) ellipse ({sin{20}} and {sin{20}*0.2});

                    \draw[red,->,rotate=12] (0,0) -- (0,1.12) node[pos=1.1,black] {\(\hat{\bm{\alpha}}\)};

                    \draw[black,dashed,->] ([shift=(90:0.8)]0,0) arc (90:70:0.8)  node[pos=1.5,black] {\(\delta_\text{max}\)};

                    \draw[black,dashed,->] ([shift=(90:0.7)]0,0) arc (90:102:0.7)  node[pos=1.3,black] {\(\delta\)};

                    \draw[black,->] (0,\armLength) -- (0,0) node[pos=1.12,black] {\(\bm{r}_T\)};
                \end{scope} 
            \end{scope}
        \end{scope}

        \draw[black,->] (0,0) -- (\xB*0.98,\yB*0.98) node[above,pos=0.5,black] {\(\bm{r}\)};

        \draw[black,dashed,->] ([shift=(180:1.8)]0,0) arc (180:142:1.8)  node[pos=1.2,black] {\(\gamma\)};
        
    \end{tikzpicture}
    \caption{Definitions of the inertial and body  reference frames and some of the variables and parameters.}
    \vspace{-5mm}
    \label{fig: frames schematic}
\end{figure}
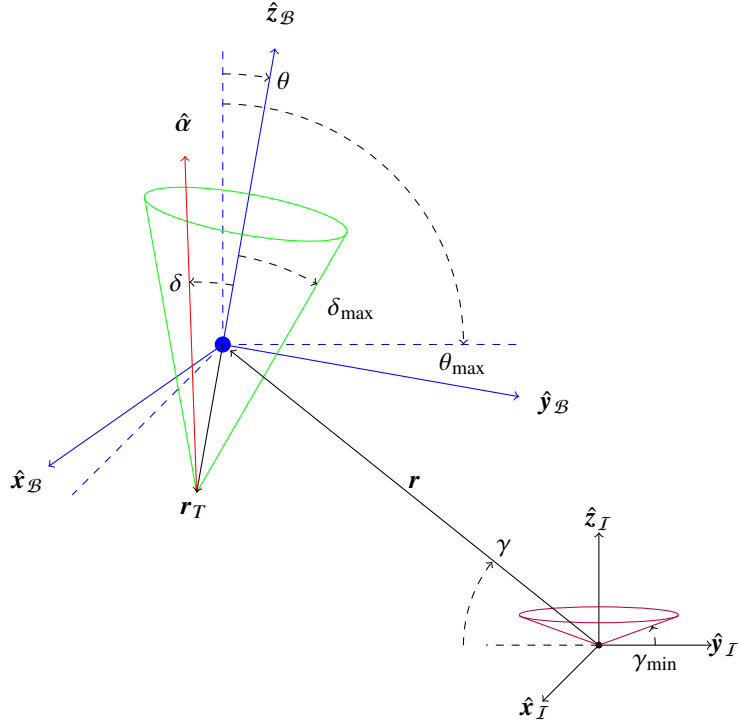

The dynamical model and accompanying constraints are depicted in the schematic in Fig.~\ref{fig: frames schematic}. The translational dynamics of the rocket are modeled in an inertial reference frame, \(\mathcal{F}_\mathcal{I}=\left\{\hat{\bm{x}}_\mathcal{I},\hat{\bm{y}}_\mathcal{I},\hat{\bm{z}}_\mathcal{I}\right\}\), centered on the surface of a flat and non-rotating Earth. We also define a body reference frame centered at the rocket's center of mass, \(\mathcal{F}_\mathcal{B}=\left\{\hat{\bm{x}}_\mathcal{B},\hat{\bm{y}}_\mathcal{B},\hat{\bm{z}}_\mathcal{B}\right\}\), where the \(\hat{\bm{z}}_\mathcal{B}\) axis points along the nose of the rocket. The rocket's center of mass location and thus body frame origin location is assumed to be constant relative to the rocket. The state of the rocket's body frame origin with respect to the origin of the inertial frame is defined by Cartesian position and velocity vectors denoted by \(\bm{r}=\left[r_x,r_y,r_z\right]^\top\) and \(\bm{v}=\left[v_x,v_y,v_z\right]^\top\), respectively. The rocket is subject to the constant acceleration of gravity denoted by \(\bm{g}=\left[0,0,-g_0\right]^\top\) where \(g_0\) denotes the sea-level gravitational acceleration of the Earth. The rocket is assumed to have a single thruster that can gimbal. The gimbal point is located by the position vector \(\bm{r}_T=\left[r_{T,x},r_{T,y},r_{T,z}\right]^\top\), expressed in the body frame. The thrust magnitude, $T$, is constrained as \(T\in\left[T_\text{min},T_\text{max}\right]\), where \(T_\text{min}\) and \(T_\text{max}\) denote the minimum and maximum thrust magnitudes, respectively. The thrust steering direction is denoted by the unit vector \(\hat{\bm{\alpha}}=\left[\alpha_x,\alpha_y,\alpha_z\right]^\top\), which is expressed in the body frame. This particular magnitude-vector control parametrization is different from Refs. \cite{szmuk_successive_2016,sagliano_six-degree--freedom_2024} and makes derivation of the necessary conditions of optimality more straightforward.

Under rigid-body kinematics assumptions, the orientation of the body frame is modeled with unit quaternions, \(\bm{q}=\left[q_0,q_1,q_2,q_3\right]^\top\). The angular velocity of the body frame relative to the inertial frame is expressed in the body frame as \(\bm{\omega}=\left[\omega_x,\omega_y,\omega_z\right]^\top\). The moment of inertia matrix is denoted by \(\bm{J}=\diag{\left(\left[J_x,J_y,J_z\right]\right)}\). The transformation from the inertial to the body frame is defined by the direction cosine matrix (DCM) in terms of quaternions as,
\begin{equation}
    \bm{C}_{\mathcal{B}\leftarrow\mathcal{I}} = \begin{bmatrix} 1 - 2(q_2^2 + q_3^2) & 2(q_1q_2+q_0q_3) & 2(q_1q_3 - q_0q_2) \\ 2(q_1q_2-q_0q_3) & 1-2(q_1^2+q_3^2) & 2(q_2q_3+q_0q_1) \\ 2(q_1q_3+q_0q_2) & 2(q_2q_3-q_0q_1) & 1-2(q_1^2+q_2^2) \end{bmatrix},
\end{equation}
and its transpose, \(\bm{C}_{\mathcal{B}\leftarrow\mathcal{I}}^\top=\bm{C}_{\mathcal{I}\leftarrow\mathcal{B}}\), denotes transformation from the body frame to the inertial frame. Following Ref. \cite{sagliano_six-degree--freedom_2024}, a drag force, \(\bm{D}\), is considered which is defined with respect to the inertial frame as $\bm{D} = -\frac{1}{2}\rho\left\|\bm{v}\right\|\bm{v}SC_D$, where \(\rho\) denotes a constant atmospheric density, \(S\) denotes the aerodynamic reference area, and \(C_D\) denotes a constant drag coefficient. The drag force will only affect the translational dynamics. Aerodynamic moments can be included in our methodology without loss of generality, similar to how they are modeled in Ref.~\cite{szmuk_successive_2018}. The rocket's mass is denoted by \(m\) and the specific impulse of the main engine is denoted by \(I_\text{sp}\). The state vector for the rocket is denoted by \(\bm{x}=\left[\bm{r}^\top,\bm{v}^\top,\bm{q}^\top,\bm{\omega}^\top,m\right]^\top \in \mathbb{R}^{14}\). The 6DOF equations of motion can be written as,
\begin{align}\label{eq: eom}
    \dot{\bm{r}} & = \bm{v}, & 
    \dot{\bm{v}} & = \bm{g} + \frac{T}{m}\bm{C}_{\mathcal{I}\leftarrow\mathcal{B}}\hat{\bm{\alpha}} + \frac{\bm{D}}{m}, &
    \dot{\bm{q}} & = \frac{1}{2}\bm{\Omega}\bm{q}, \nonumber \\ 
    \dot{\bm{\omega}} & = \bm{J}^{-1}\left(\bm{r}_T^\times T\hat{\bm{\alpha}} - \bm{\omega}^\times \bm{J}\bm{\omega}\right), & 
    \dot{m} & = -\frac{T}{I_\text{sp}g_0},
\end{align}
where superscript `\(\times\)' denotes the skew-symmetric matrix forms of $\bm{\omega}$ and $\bm{r}_T$ and \(\bm{\Omega}\) denotes a \(4\times4\) skew-symmetric matrix of \(\bm{\omega}\) written as,
\begin{align}
    \bm{\omega}^\times & = \begin{bmatrix} 0 & -\omega_z & \omega_y \\ \omega_z & 0 & -\omega_x \\ -\omega_y & \omega_x & 0 \end{bmatrix}, & \bm{r}_T^\times & = \begin{bmatrix} 0 & -r_{T,z} & r_{T,y} \\ r_{T,z} & 0 & -r_{T,x} \\ -r_{T,y} & r_{T,x} & 0 \end{bmatrix}, &  \bm{\Omega}=\begin{bmatrix} 0 & -\omega_x & -\omega_y & -\omega_z \\ \omega_x & 0 & \omega_z & -\omega_y \\ \omega_y & -\omega_z & 0 & \omega_x \\ \omega_z & \omega_y & -\omega_x & 0 \end{bmatrix}.
\end{align}

\section{Constrained Fuel- and Time-Optimal Trajectory Optimization Problems}\label{sec: optimal control problem}

The OCP is to find the thrust magnitude, \(T\), and steering control, \(\hat{\bm{\alpha}}\), time histories that land the rocket upright starting from some initial state with a free time-of-flight, \(t_f\), in either minimal propellant mass (i.e., fuel-optimal) or in minimal time-of-flight (i.e., time-optimal). The problem is solved over the time horizon \(t\in[0,t_f]\). The initial and final states, \(\bm{x}(0)\) and \(\bm{x}(t_f)\), are summarized as,
\begin{align} 
    \bm{r}(0) & = \bm{r}_0, &
    \bm{v}(0) & = \bm{v}_0, &
    \bm{q}(0) & = \text{free}, &
    \bm{\omega}(0) & = \bm{\omega}_{0}, &
    m(0) & = m_0, \label{eq: initial BCs} \\
    \bm{r}(t_f) & = \bm{r}_f, &
    \bm{v}(t_f) & = \bm{v}_f, &
    \bm{q}(t_f) & = \bm{q}_{f} & 
    \bm{\omega}(t_f) & = \bm{\omega}_{f}, &
    m(t_f) & = \text{free}, \label{eq: final BCs}
\end{align}
where the initial orientation and final mass are free. The rocket is subject to the dynamical constraints defined in Eq.~\eqref{eq: eom}. The angular velocity magnitude, \(\left\|\bm{\omega}\right\|\), and the glideslope angle, \(\gamma\), are constrained as, 
\begin{align}
    \left\|\bm{\omega}\right\| & \leq \omega_\text{max}, \label{eq: angular velocity constraint} \\  \gamma_\text{min} \leq \gamma & =  \tan^{-1}{\left(\frac{r_z}{\sqrt{r_x^2+r_y^2}}\right)},\label{eq: glideslope constraint}
\end{align}
where \(\omega_\text{max}\) and $\gamma_\text{min}$ denote the bounds. The rocket's tilt angle, \(\theta\) (i.e., the angle between \(\hat{\bm{z}}_\mathcal{B}\) and \(\hat{\bm{z}}_\mathcal{I}\) as it is shown in Fig.~\ref{fig: frames schematic}), is constrained to be less than some maximum tilt angle, \(\theta_\text{max}\). Similar to Refs. \cite{szmuk_successive_2016,sagliano_six-degree--freedom_2024} we define this constraint using the dot product operation performed in the inertial coordinate frame, i.e.,
\begin{equation}\label{eq: tilt angle constraint}
    \theta = \cos^{-1}{\left(\hat{\bm{z}}_\mathcal{I}^\top \bm{C}_{\mathcal{I}\leftarrow\mathcal{B}} \hat{\bm{z}}_\mathcal{B}\right)} = \cos^{-1}{\left(1-2\left(q_1^2+q_2^2\right)\right)} \leq \theta_\text{max}.
\end{equation}

The thrust magnitude, \(T\), steering vector, $\hat{\bm{\alpha}}$, and thruster's gimbal angle, $\delta$ (see Fig.~\ref{fig: frames schematic}), are constrained as, 
\begin{align} \label{eq: control constraints}
    T_\text{min} & \leq T \leq T_\text{max}, & \hat{\bm{\alpha}}^\top\hat{\bm{\alpha}} & = 1, &   \delta = \cos^{-1}{\left(\alpha_z\right)} \leq \delta_\text{max},
\end{align}
with \(\delta_\text{max}\), denoting a maximum gimbal angle value. The fuel- and time-optimal trajectory optimization problems are summarized in Eq.~\eqref{eq: ocp} as,

\begin{align}\label{eq: ocp}
    & \mathcal{P}_\text{FO} \left\{
    \begin{aligned}
        & \min_{T,\hat{\bm{\alpha}},t_f} \quad -m(t_f) \\
        & \text{s.t.,} \quad \text{Eqs.~}\eqref{eq: eom},\eqref{eq: initial BCs},\eqref{eq: final BCs},\eqref{eq: angular velocity constraint},\eqref{eq: glideslope constraint},\eqref{eq: tilt angle constraint},\eqref{eq: control constraints}
    \end{aligned} \right.
    & \mathcal{P}_\text{TO} \left\{
    \begin{aligned}
        & \min_{T,\hat{\bm{\alpha}},t_f} \quad t_f \\
        & \text{s.t.,} \quad \text{Eqs.~}\eqref{eq: eom},\eqref{eq: initial BCs},\eqref{eq: final BCs},\eqref{eq: angular velocity constraint},\eqref{eq: glideslope constraint},\eqref{eq: tilt angle constraint},\eqref{eq: control constraints}
    \end{aligned} \right. 
\end{align}

\section{Derivation of Necessary Conditions}\label{sec: derivation of necessary conditions}

The indirect method is used to find stationary solutions to the OCPs in Eq.~\eqref{eq: ocp} through the first-order necessary conditions of optimality. Following the notation of Ref. \cite{bryson_applied_1975}, the SOPICs in Eqs.~\eqref{eq: angular velocity constraint}, \eqref{eq: glideslope constraint}, and \eqref{eq: tilt angle constraint} are rewritten in the form \(S_i\leq0\),
\begin{align} \label{eq: new control constraints}
    S_1 & = \bm{\omega}^\top\bm{\omega} - \omega_\text{max}^2 \leq 0, &
    S_2 & = \gamma_\text{min} - \gamma \leq 0, & 
    S_3 & = \theta - \theta_\text{max} \leq 0. 
\end{align}

Note that the angular velocity magnitude constraint is rewritten without the square root for a better regularization of its first derivative with respect to states, which is a step in the derivation of the costate (adjoint) equations. 

A penalty function approach is taken to enforce the SOPICs. Penalty function approaches offer promising results for solving practical/challenging OCPs \cite{epenoy_fuel_2011,mall_uniform_2020}. We use secant penalty functions and augment them to the Lagrangian, \(\mathcal{L}\). The secant penalty function is employed to enforce the \(i\)-th SOPIC, \(S_i\leq 0\), as \(\rho_i\cdot\sec{\left(\pi/2\cdot P_i\right)}\) where \(P_i\) is \(S_i\) rewritten in the form \(P_i\leq 1\). Despite the secant function being singular at \(P_i=1\), the parameter \(\rho_i\) is introduced to approximate the active constraint arcs by starting with a high parameter value and then gradually lowering it to near 0 through continuation, i.e., \(\rho_i\ \rightarrow 0\). This regularization has worked successfully, i.e., \(S_i\rightarrow0\) as \(\rho_i \rightarrow 0\) on the active constraint arcs, but we don't rigorously show this (see Refs. \cite{epenoy_fuel_2011} and \cite{malisani_interior_2016} for convergence proofs of interior penalty methods). 

Let the secant penalty functions that are augmented to the Lagrangian be defined as,
\begin{align} \label{eq: secant penalty functions}
    \tilde{S}_1 & = \sec{\left(\frac{\pi}{2}\cdot\frac{\bm{\omega}^\top\bm{\omega}}{\omega_\text{max}^2}\right)}, & \tilde{S}_2 & = \sec{\left(\frac{\pi}{2} \cdot \frac{\gamma_\text{min}}{\gamma}\right)}, & \tilde{S}_3 & = \sec{\left(\frac{\pi}{2}\cdot \frac{\theta}{\theta_\text{max}}\right)}.
\end{align}

Let \(\rho_\omega\), \(\rho_\gamma\), and \(\rho_\theta\) denote weighting coefficients associated with the functions, respectively, in Eq.~\eqref{eq: secant penalty functions}. The Lagrangian is then defined as,
\begin{align}
    \mathcal{L} = \rho_\omega \tilde{S}_1 + \rho_\gamma \tilde{S}_2 + \rho_\theta \tilde{S}_3.
\end{align}

The weighting coefficients will be used within a numerical continuation scheme to solve the resulting boundary-value problems, as explained later. Let the costate vector be denoted by \(\bm{\lambda}=\left[\bm{\lambda}_{\bm{r}}^\top,\bm{\lambda}_{\bm{v}}^\top,\bm{\lambda}_{\bm{q}}^\top,\bm{\lambda}_{\bm{\omega}}^\top,\lambda_m\right]^\top \in \mathbb{R}^{14}\). The (optimal control/variational) Hamiltonian can be written as,
\begin{equation} \label{eq: hamiltonian}
    H = \mathcal{L} + \bm{\lambda}_{\bm{r}}^\top\bm{v} + \bm{\lambda}_{\bm{v}}^\top\left(\bm{g} + \frac{T}{m}\bm{C}_{I\leftarrow B}\hat{\bm{\alpha}} + \frac{\bm{D}}{m}\right) + \frac{1}{2}\bm{\lambda}_{\bm{q}}^\top\bm{\Omega}\bm{q} + \bm{\lambda}_{\bm{\omega}}^\top\bm{J}^{-1}\left(\bm{r}_T^\times T\hat{\bm{\alpha}} - \bm{\omega}^\times \bm{J}\bm{\omega}\right) - \lambda_m\frac{T}{I_\text{sp}g_0}.
\end{equation}

The costate equations can be derived using the Euler-Lagrange equation, $\dot{\bm{\lambda}} = -(\partial H/\partial \bm{x})^\top$, which is constructed with automatic differentiation (AD) using CasADi \cite{andersson_casadi_2019}. Based on our experiences \cite{nurre_end--end_2024}, using CasADi significantly facilitates the implementation of indirect methods. Note that the derivative of the tilt angle, \(\theta\) (see Eq.~\eqref{eq: tilt angle constraint}), with respect to both \(q_1\) and \(q_2\) (which occurs in $\dot{\bm{\lambda}} = -(\partial H/\partial \bm{x})^\top$) has a singularity when \(q_1=q_2=0\) i.e., 
\begin{align}
    \frac{\partial \theta}{\partial q_1} & = \frac{4q_1}{\sqrt{1-(2(q_1^2+q_2^2)-1)^2}}, & 
    \frac{\partial \theta}{\partial q_2} & = \frac{4q_2}{\sqrt{1-(2(q_1^2+q_2^2)-1)^2}}.
\end{align}

To avoid this, we slightly modify the final orientation boundary condition in Table \ref{tab: boundary conditions} with no impact on the solutions.

\begin{remark}
The Lagrange multipliers associated with the SOPICs following the direct adjoining approach \cite{hartl_survey_1995,ross_primer_2015} can be calculated a posteriori under the penalty function approach and verified to satisfy the complementarity conditions. We next show how we arrive at the expression for these ``approximate Lagrange multipliers,'' denoted by \(\tilde{\eta}\). Because we only empirically verify these Lagrange multipliers (as will be shown in Section \ref{sec: results}) and do not provide a rigorous proof of equivalence, we present it as a conjecture.
\end{remark}

\begin{conjecture} \label{conjecture 1}
    Let \(\eta_i\) be the Lagrange multiplier associated with the \(i\)-th SOPIC, \(S_i\leq 0\). Following the direct adjoining approach to enforcing SOPICs (see Informal Theorem 4.1 in Ref. \cite{hartl_survey_1995}), \(S_i\) is adjoined to the Hamiltonian with \(\eta_i\). An accurate approximation of \(\eta_i\), which is denoted as \(\tilde{\eta}_i\), can be derived by equating the Hamiltonian of the direct adjoining approach to the Hamiltonian of the secant-function-augmented Hamiltonian, i.e.,
    \begin{equation} \label{eq: eta tilde}
        \bm{\lambda}^\top\bm{f} + \eta_iS_i \approx \bm{\lambda}^\top\bm{f}  +\rho_i\sec{\left(\frac{\pi}{2} P_i\right)}, \rightarrow \tilde{\eta}_i = -\rho_i\sec{\left(\frac{\pi}{2} P_i\right)}/S_i,
    \end{equation}
    where \(P_i\) is the constraint \(S_i\leq0\) redefined in the form \(P_i\leq1\) and the approximately-equal-sign is used since we don't specify \(\rho_i\) here. Also, the negative sign is included since we express the constraint in the form \(S_i\leq0\), leaving \(S_i\) always nonpositive and causing \(\tilde{\eta}_i\) to otherwise violate the complementarity condition under the minimum principle. After a solution is found, \(\tilde{\eta}_i\) can be calculated using Eq.~\eqref{eq: eta tilde} and verified that it satisfies complementarity slackness conditions. To empirically show that \(\tilde{\eta}_i \rightarrow \eta_i\) as \(\rho_i \rightarrow 0\), we use DIDO \cite{ross_enhancements_2020} to solve the OCPs and retrieve \(\eta_i\) from its solution.
\end{conjecture}

The control equality/inequality constraints in Eqs.~\eqref{eq: control constraints} are enforced by directly adjoining them to the Hamiltonian with additional Lagrange multipliers to form the augmented Hamiltonian, \(\bar{H}\). This allows us to algebraically solve for exact open-loop control expressions (i.e., a function of states and costates) using the strong form of optimality and the complementarity conditions \cite{bryson_optimal_1963}. Eqs.~\eqref{eq: control constraints} are rewritten in the form \(C\leq0\) as,
\begin{align} 
    C_1 & = T_\text{min} - T \leq 0,  & 
    C_2 & = T - T_\text{max} \leq 0,  &
    C_3 & = \hat{\bm{\alpha}}^\top\hat{\bm{\alpha}} - 1 = 0, &
    C_4 & = \cos{\left(\delta_\text{max}\right)} - \hat{\bm{z}}_\mathcal{B}^\top\hat{\bm{\alpha}} \leq 0. 
\end{align}

Let \(\mu_1\), \(\mu_2\), \(\mu_3\), and \(\mu_4\) be the Lagrange multipliers associated with $C_1$, $C_2$, $C_3$, and $C_4$, respectively. The augmented Hamiltonian is then expressed as,
\begin{equation} \label{eq: augmented hamiltonian}
\bar{H} = H + \mu_1C_1 + \mu_2C_2 + \mu_3C_3 + \mu_4C_4.
\end{equation}

Closed-form expressions for the control variables \(T\) and \(\hat{\bm{\alpha}}\) will now be derived. Eq.~\eqref{eq: augmented hamiltonian} is rewritten with all control-independent terms collected into \(H_0 = H_0(\bm{x},\bm{\lambda})\) as,
\begin{equation}
    \bar{H} = H_0 + \bm{\lambda}_{\bm{v}}^\top\frac{T}{m}\bm{C}_{I\leftarrow B}\hat{\bm{\alpha}} + \bm{\lambda}_{\bm{\omega}}^\top\bm{J}^{-1}\bm{r}_T^\times T\hat{\bm{\alpha}} - \lambda_m\frac{T}{I_\text{sp}g_0} - \mu_1T + \mu_2T + \mu_3\hat{\bm{\alpha}}^\top\hat{\bm{\alpha}} - \mu_4\hat{\bm{z}}_\mathcal{B}^\top\hat{\bm{\alpha}}.
\end{equation}

Let \(\bm{p}=\left[p_x,p_y,p_z\right]^\top\) be defined as
\begin{equation} \label{eq: primer vector}
    \bm{p}^\top = \frac{\bm{\lambda}_{\bm{v}}^\top\bm{C}_{I\leftarrow B}}{m} + \bm{\lambda}_{\bm{\omega}}^\top\bm{J}^{-1}\bm{r}_T^\times,
\end{equation}
which is the direction opposite of the primer vector \cite{russell_primer_2007}, i.e., the unconstrained optimal direction of thrusting.

Using Eq.~\eqref{eq: primer vector}, the augmented Hamiltonian is rewritten as,
\begin{equation} \label{eq: reduced augmented hamiltonian}
    \bar{H} = H_0 + T\bm{p}^\top\hat{\bm{\alpha}}  - \lambda_m\frac{T}{I_\text{sp}g_0} - \mu_1T + \mu_2T + \mu_3\hat{\bm{\alpha}}^\top\hat{\bm{\alpha}} - \mu_4\hat{\bm{z}}_\mathcal{B}^\top\hat{\bm{\alpha}}.
\end{equation}

The thrust magnitude control, \(T\), is derived by applying the strong form of optimality to Eq.~\eqref{eq: reduced augmented hamiltonian},
\begin{align}
    \frac{\partial \bar{H}}{\partial T} = \bm{p}^\top\hat{\bm{\alpha}} - \frac{\lambda_m}{I_\text{sp}g_0} - \mu_1 + \mu_2 = 0.
\end{align}

Due to the complementarity conditions (i.e., \(\mu_i C_i = 0 \wedge \mu_i\geq 0\) for \(i=1,2\)), only one of the three following conditions is true at any time along an extremal solution: 1) \(\mu_1 = 0\) and \(\mu_2\neq0\), 2) \(\mu_1 \neq 0\) and \(\mu_2 = 0\), or 3) \(\mu_1 = \mu_2 = 0\). Thus, \(T\) will be bang-bang (conditions 1 and 2) unless a singular arc occurs (condition 3). It can be verified that the optimal thrust magnitude control, \(T^*\), can be written as,
\begin{align} \label{eq: optimal thrust magnitude}
    T^* & \begin{cases} = T_\text{max} & S_T > 0, \\ \in [T_\text{min},T_\text{max}] & S_T = 0, \\ = T_\text{min} & S_T < 0,
    \end{cases} & S_T & = -\bm{p}^\top\hat{\bm{\alpha}} + \frac{\lambda_m}{I_\text{sp}g_0},
\end{align}
where \(S_T\) denotes the so-called  thrust magnitude switching function. 

The optimal expression for \(\hat{\bm{\alpha}}\) is not known at this point which  appears in \(S_T\) in Eq.~\eqref{eq: optimal thrust magnitude}. However, its closed-form expression, as will be shown, is independent of \(T^*\) and can be substituted into \(S_T\). To derive the optimal steering control, $\hat{\bm{\alpha}}$, the strong form of optimality is again invoked on Eq.~\eqref{eq: reduced augmented hamiltonian} as,
\begin{equation} \label{eq: mu3}
    \frac{\partial \bar{H}}{\partial \hat{\bm{\alpha}}} = T\bm{p}^\top + \mu_3\hat{\bm{\alpha}}^\top - \mu_4\hat{\bm{z}}_\mathcal{B}^\top = \bm{0} \rightarrow \hat{\bm{\alpha}} = -\frac{T\bm{p} - \mu_4\hat{\bm{z}}_\mathcal{B}}{\mu_3}, \xrightarrow[]{\text{substitute into}~C_3=0} \mu_3 = \pm \left\|T\bm{p} - \mu_4\hat{\bm{z}}_\mathcal{B}\right\|.
\end{equation}

To determine that the positive sign in Eq.~\eqref{eq: mu3} should be taken, the Legendre-Clebsch condition can be invoked on Eq.~\eqref{eq: reduced augmented hamiltonian}, \(\partial^2 \bar{H}/\partial \hat{\bm{\alpha}}^2 = \mu_3 \geq 0\). The thrust steering control expression can then be written as,
\begin{equation} \label{eq: alpha with mu4}
    \hat{\bm{\alpha}} = -\frac{T\bm{p} - \mu_4\hat{\bm{z}}_\mathcal{B}}{\left\|T\bm{p} - \mu_4\hat{\bm{z}}_\mathcal{B}\right\|}.
\end{equation}

\begin{figure}[]
    \centering
    \begin{tikzpicture}
        \begin{scope}[scale=4]                    
                \draw[blue,->] (0,0) -- (1,0) node[pos=1.12,black] {\(\hat{\bm{y}}_\mathcal{B}\)};
                \draw[blue,->] (0,0) -- (0,0) -- (0,1) node[pos=1.12,black] {\(\hat{\bm{z}}_\mathcal{B}\)};
                \draw[blue,->] (0,0) -- (0,0) -- (-0.5,-0.5) node[pos=1.12,black] {\(\hat{\bm{x}}_\mathcal{B}\)};

                \draw[black,dashed,->] ([shift=(90:0.6)]0,0) arc (90:110:0.6)  node[pos=1.5,black] {\(\delta_\text{max}\)};

                \filldraw[blue] (0,0) circle (0.025);
                
                \draw[green,rotate=20] (0,0) -- (0,0.8);
                \draw[green,rotate=-20] (0,0) -- (0,0.8);
                \draw[green] (0,0.75) ellipse ({0.8*sin{20}} and {0.8*sin{20}*0.2});

                \draw[red,->] (0,0) -- (sin{20},cos{20}) node[pos=1.1,black] {\(\hat{\bm{\alpha}}^*(\mu_4^+)\)};
                \draw[red,->] (0,0) -- (sin{20},-cos{20}) node[pos=1.1,black] {\(\hat{\bm{\alpha}}^*(\mu_4^-)\)};
                \draw[red,dashed] (sin{20},cos{20}) -- (sin{20},-cos{20}) node[pos=0.85,black,right] {\(-\cot{\left(\delta_\text{max}\right)}\sqrt{p_x^2+p_y^2}\)} node[pos=0.35,black,right] {\(\cot{\left(\delta_\text{max}\right)}\sqrt{p_x^2+p_y^2}\)} ;
                \filldraw[red] (sin{20},-0.237) circle (0.02);

                \draw[red,->] (0,0) -- (sin{40},cos{40}) node[pos=1.1,black] {\(-\frac{\bm{p}}{\left\|\bm{p}\right\|}\)};
                \draw[red,dashed] (sin{40},cos{40}) -- (sin{40},-cos{40}*0.6);
                \draw[red,dashed] (sin{40},-cos{40}*0.6) -- (0,0);
                
        \end{scope}
    \end{tikzpicture}
    \caption{Schematic of the gimbal constraint.}
    \label{fig: gimbal constraint}
\end{figure}
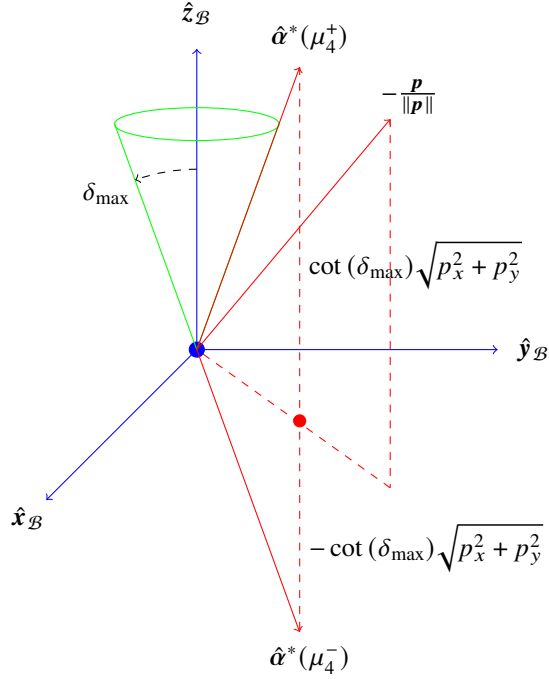

The variable \(\mu_4\) remains to be determined. When the gimbal constraint ($C_4$) is not active, then, according to complementarity, \(\mu_4=0\) and the optimal steering expression is simply in the direction of the primer vector, i.e., \(-\bm{p}/\left\|\bm{p}\right\|\), which is consistent with other unconstrained steering problems in the literature. An expression for \(\mu_4\), however, must be derived when the gimbal constraint is active and $C_4 = 0$. The expression resulting from substituting Eq.~\eqref{eq: alpha with mu4} into $C_4 = 0$, becomes
\begin{equation} \label{eq: c4 with alpha and mu4}
    \cos{\left(\delta_\text{max}\right)} + \hat{\bm{z}}_\mathcal{B}^\top\frac{T\bm{p} - \mu_4\hat{\bm{z}}_\mathcal{B}}{\left\|T\bm{p} - \mu_4\hat{\bm{z}}_\mathcal{B}\right\|}=0,
\end{equation}
and is quadratic in \(\mu_4\) and has two solutions written as,
\begin{align} \label{eq: mu4 roots}
    \mu_4^+ & = p_z + \cot{\left(\delta_\text{max}\right)\sqrt{p_x^2+p_y^2}}, & \mu_4^- & = p_z - \cot{\left(\delta_\text{max}\right)\sqrt{p_x^2+p_y^2}},
\end{align}
giving two possible steering control expressions for when the gimbal constraint is active, i.e.,
\begin{align} \label{eq: steering with both mu4}
    \frac{\hat{\bm{\alpha}}(\mu_4^+)}{\|\cdot\|} & = \begin{bmatrix} -p_x ,-p_y,  \cot{\left(\delta_\text{max}\right)\sqrt{p_x^2+p_y^2}} \end{bmatrix}, & \frac{\hat{\bm{\alpha}}(\mu_4^-)}{\|\cdot\|} & = \begin{bmatrix} -p_x,-p_y,-\cot{\left(\delta_\text{max}\right)\sqrt{p_x^2+p_y^2}} \end{bmatrix},
\end{align}
where \(\|\cdot\|\) is the 2-norm of the expression in the denominators of Eqs.~\eqref{eq: c4 with alpha and mu4} and \eqref{eq: steering with both mu4}. 
To choose the correct expression, we can consider Figure \ref{fig: gimbal constraint}; when the gimbal constraint is active, \(\mu_4^+\) corresponds to  the thrust steering direction lying on the constraint boundary, whereas \(\mu_4^-\) will be the same vector but flipped about the \(\hat{\bm{x}}_\mathcal{B}\)-\(\hat{\bm{y}}_\mathcal{B}\) axis. The optimal value for \(\mu_4\), \(\mu_4^*\), can therefore be defined in a piecewise manner to have a closed-form expression for \(\hat{\bm{\alpha}}\). Let \(S_\delta\) denote a so-called gimbal constraint-activation switching function, \(\mu_4^*\) can then be defined as,
\begin{align} \label{eq: optimal mu4}
    \mu_4^* &  \begin{cases} = 0 & S_\delta > 0, \\ = \mu_4^+ & S_\delta \leq 0,
    \end{cases} & S_\delta & = p_z - \cot{\left(\delta_\text{max}\right)}\sqrt{p_x^2+p_y^2}.
\end{align}

Upon substituting Eq.~\eqref{eq: optimal mu4} into Eq.~\eqref{eq: alpha with mu4}, the optimal steering expression can be derived as,
\begin{equation} \label{eq: optimal steering}
\hat{\bm{\alpha}}^* = -\frac{\bm{p}-\mu_4^*\hat{\bm{z}}_\mathcal{B}}{\left\|\bm{p}-\mu_4^*\hat{\bm{z}}_\mathcal{B}\right\|}.
\end{equation}
Note that this expression given in Eq.~\eqref{eq: optimal steering} is independent of \(T\), allowing us to substitute it back into Eq.~\eqref{eq: optimal thrust magnitude}. 

The remaining necessary conditions of optimality are the transversality conditions due to the free terminal states in Eqs.~\eqref{eq: initial BCs} and \eqref{eq: final BCs} and the stationarity condition resulting from the fact that the problem is free-final time. Applying the transversality conditions to the free initial orientation we have $\bm{\lambda}_{\bm{q}}(0) = \bm{0}$. For the fuel-optimal problem, the mass transversality condition and free-final-time stationarity conditions are $\lambda_m(t_f) = -1$ and $H(t_f) = 0$, whereas, for the time-optimal problem we have $\lambda_m(t_f) = 0$ and $H(t_f) = -1$. The Hamiltonian is not an explicit function of time and its value will remain constant. This property can be used to verify for the optimality \textit{a posteriori}. Note that we treat the final time, \(t_f\), as a parameter and use the necessary conditions from Chapter 14 in Ref. \cite{hull_optimal_2003}. If the Hamiltonian does not remain constant over time, then the resulting Hamiltonian condition would be different.

\subsection{Regularization of Control Expressions}

The piecewise-continuous expressions \(T^*\) (Eq.~\eqref{eq: optimal thrust magnitude}) and \(\mu_4^*\) (Eq.~\eqref{eq: optimal mu4}) reduce the domain of convergence of the resulting boundary-value problems \cite{bertrand_new_2002,petukhov_method_2012}. Thus, they are embedded into one-parameter families of smooth curves using an L2 norm-based regularization \cite{taheri_l2_2023}. The new regularized expressions, denoted using \(\tilde{\cdot}\), are defined as,
\begin{equation} \label{eq: regularized thrust}
    \tilde{T}^* = \frac{1}{2} \left [ \left(T_\text{max} + T_\text{min}\right) + \left(T_\text{max} - T_\text{min}\right)\frac{S_T}{\sqrt{S_T^2+\rho_T^2}} \right ],
\end{equation}
\begin{equation} \label{eq: regularized steering}
    \tilde{\mu}_4^* = \frac{\mu_4^+}{2}\left(1 + \frac{S_\delta}{\sqrt{S_\delta^2+\rho_\delta^2}}\right) \implies \tilde{\hat{\bm{\alpha}}}^* = -\frac{\bm{p}-\tilde{\mu}_4^*\hat{\bm{z}}_\mathcal{B}}{\left\|\bm{p}-\tilde{\mu}_4^*\hat{\bm{z}}_\mathcal{B}\right\|}.
\end{equation}

The continuation parameters \(\rho_T\) and \(\rho_\delta\) are introduced such that \(\tilde{T}^*\) and \(\tilde{\mu}_4^*\) will converge to \(T^*\) and \(\mu_4^*\) as \(\rho_T\rightarrow0\) and \(\rho_\delta\rightarrow0\). Moreover, the L2 norm-based regularization can approximate singular arcs in \(T^*\) \cite{taheri_l2_2023} should they arise.

\section{Fuel- and Time-Optimal Regularized Two-Point Boundary-Value Problems} \label{sec: boundary value problem}

The necessary conditions, with regularized control expressions and SOPICs implicitly enforced with penalty functions, transform what would be a MPBVP into a multi-parameter smooth family of neighboring two-point boundary-value problems. The independent variable, \(t\in[0,t_f]\), is scaled by the time of flight, \(t_f\), and redefined as \(\tau\in[0,1]\). This new system is stated in the state-space form as,
\begin{align} \label{eq: bvp system}
    \bm{z}(\tau) & = \begin{bmatrix} \bm{x}(\tau) \\ \bm{\lambda}(\tau) \end{bmatrix}, & \dot{\bm{z}}(\tau) & = \begin{bmatrix} \dot{\bm{x}}(\tau) \\ \dot{\bm{\lambda}}(\tau) \end{bmatrix} = t_f \bm{f}(\tau ,\bm{z}(\tau);t_f,\bm{\rho}),
\end{align}
where $\bm{f} = \left[\dot{\bm{r}}^\top,\dot{\bm{v}}^\top,\dot{\bm{q}}^\top,\dot{\bm{\omega}}^\top,\dot{m}\right]^\top$ and \(\bm{\rho} = \left[\rho_T,\rho_\delta,\rho_\omega,\rho_\gamma,\rho_\theta\right]^\top\). The initial conditions are 
\begin{align}
    \bm{r}(0) & = \bm{r}_0, &
    \bm{v}(0) & = \bm{v}_0, &
    \bm{q}(0) & = \text{free}, &
    \bm{\omega}(0) & = \bm{\omega}_{0}, &
    m(0) & = m_0, \\
    \bm{\lambda}_{\bm{r}}(0) & = \text{free}, &
    \bm{\lambda}_{\bm{v}}(0) & = \text{free}, &
    \bm{\lambda}_{\bm{q}}(0) & = \bm{0}, &
    \bm{\lambda}_{\bm{\omega}}(0) & = \text{free}, &
    \bm{\lambda}_m(0) & = \text{free},
\end{align}
and the final boundary conditions are written as,
\begin{align}
    \bm{r}(1) & = \bm{r}_f, &
    \bm{v}(1) & = \bm{v}_f, &
    \bm{q}(1) & = \bm{q}_f, &
    \bm{\omega}(1) & = \bm{\omega}_f, &
    m(1) & = \text{free}, \\
    \bm{\lambda}_{\bm{r}}(1) & = \text{free}, &
    \bm{\lambda}_{\bm{v}}(1) & = \text{free}, &
    \bm{\lambda}_{\bm{q}}(1) & = \text{free}, &
    \bm{\lambda}_{\bm{\omega}}(1) & = \text{free}, &
    \bm{\lambda}_m(1) & = -1,
\end{align}

Formulated as a nonlinear shooting problem, there are 15 unknown variables: 14 states/costates at initial time and the time of flight. The decision vector can be written compactly as \( \bm{\Gamma} = [\bm{\lambda}_{\bm{r}}^\top(0)\), \(\bm{\lambda}_{\bm{v}}^\top(0)\), \(\bm{q}^\top(0)\), \(\bm{\lambda}_{\bm{\omega}}^\top(0)\), \(\lambda_m(0), t_f]^\top\). The number of unknowns has to be equal to the number of constraints to have a well-posed boundary-value problem. The shooting problems associated with the fuel- and time-optimal problems can be written as,
\begin{align}
    \bm{\psi}_\text{FO}(\bm{\Gamma}; \bm{\rho}) & = \begin{bmatrix} \bm{r}^\top(1) - \bm{r}^\top_f, \bm{v}^\top(1) - \bm{v}^\top_f, \bm{q}^\top(1) - \bm{q}^\top_f, \bm{\omega}^\top(1) - \bm{\omega}^\top_f, \bm{\lambda}_m(1) + 1, H(1) \end{bmatrix}^\top = \bm{0}_{15 \times 1}, \\
    \bm{\psi}_\text{TO}(\bm{\Gamma}; \bm{\rho}) & = \begin{bmatrix} \bm{r}^\top(1) - \bm{r}^\top_f, \bm{v}^\top(1) - \bm{v}^\top_f, \bm{q}^\top(1) - \bm{q}^\top_f, \bm{\omega}^\top(1) - \bm{\omega}^\top_f, \bm{\lambda}_m(1), H(1) + 1 \end{bmatrix}^\top = \bm{0}_{15 \times 1}.
\end{align}

The resulting nonlinear shooting problems don't have analytic solutions, and, instead, have to be solved numerically. This involves solving an initial-value problem (IVP) iteratively and calculating the Jacobian of the residuals with respect to the decision variables to appropriately update the unknowns.

\begin{figure}[]
    \centering
    \begin{tikzpicture}
        \begin{scope}[scale=0.9]                    
            \draw[black] (0,0) -- (4.9,0);
            \draw[black] (4.85,-0.1) -- (4.95,0.1);
            \draw[black] (5.05,-0.1) -- (5.15,0.1);
            \draw[black,->] (5.1,0) -- (11,0) node[pos=1,black,right] {\(\tau\)};
            \draw[black,->] (0,0) -- (0,5) node[pos=1,black,above] {\(\bm{z}\)};

            \draw[black] (0,-0.2) -- (0,0.2) node[pos=0,below,black] {\(\tau_0\)};
            \draw[black] (2.5,-0.2) -- (2.5,0.2) node[pos=0,below,black] {\(\tau_1\)};
            \draw[black] (7.5,-0.2) -- (7.5,0.2) node[pos=0,below,black] {\(\tau_{n-1}\)};
            \draw[black] (10,-0.2) -- (10,0.2) node[pos=0,below,black] {\(\tau_n\)};

            \def\dotRadii{0.1}
            
            \filldraw [blue] (0,1) circle (\dotRadii) node[left,black] {\(\bm{z}_1(\tau_0)\)};

            \draw[black,decoration={markings,mark=at position 0.5 with {\arrow{>}}},postaction={decorate}] (0,1) .. controls (1.5,0.25) .. (2.5,2);

            \filldraw [red] (2.5,2) circle (\dotRadii) node[below right,black] {\(\bm{z}_1(\tau_1)\)};
            \draw [decorate,decoration={brace,amplitude=4pt,mirror,raise=1ex}] (2.5,2) -- (2.5,3) node[right,pos=0.5,xshift=7pt] {\(\bm{G}_1=\bm{z}_1(\tau_1) - \bm{z}_2(\tau_1)\)};
            \filldraw [blue] (2.5,3) circle (\dotRadii) node[left,black] {\(\bm{z}_2(\tau_1)\)};

            \draw[black,decoration={markings,mark=at position 0.5 with {\arrow{>}}},postaction={decorate}] (2.5,3) .. controls (3.5, 4) .. (5,3.5);
            
            \filldraw [red] (5,3.5) circle (\dotRadii);
            \filldraw [blue] (5,4.5) circle (\dotRadii);

            \draw[black,decoration={markings,mark=at position 0.5 with {\arrow{>}}},postaction={decorate}] (5,4.5) .. controls (6, 5) .. (7.5,3.5);

            \filldraw [red] (7.5,3.5) circle (\dotRadii) node[below,black] {\(\bm{z}_{n-1}(\tau_{n-1})\)};
            \filldraw [blue] (7.5,2) circle (\dotRadii) node[below,black] {\(\bm{z}_n(\tau_{n-1})\)};

            \draw[black,decoration={markings,mark=at position 0.5 with {\arrow{>}}},postaction={decorate}] (7.5,2) .. controls (9,2) .. (10,3);

            \filldraw [red] (10,3) circle (\dotRadii) node[left,black] {\(\bm{z}_n(\tau_n)\)};
            \draw [decorate,decoration={brace,amplitude=4pt,mirror,raise=1ex}] (10,3) -- (10,4) node[right,pos=0.5,xshift=7pt] {\(\bm{G}_n=\bm{\psi}_\text{FO}~ (\text{or}~ \bm{\psi}_\text{TO})\)};
            \filldraw [blue] (10,4) circle (\dotRadii) node[left,black] {\(\bm{z}_f\)};
            
        \end{scope}
    \end{tikzpicture}
    \caption{Schematic of the indirect multiple-shooting solution scheme.}
    \label{fig: multiple shooting schematic}
\end{figure}
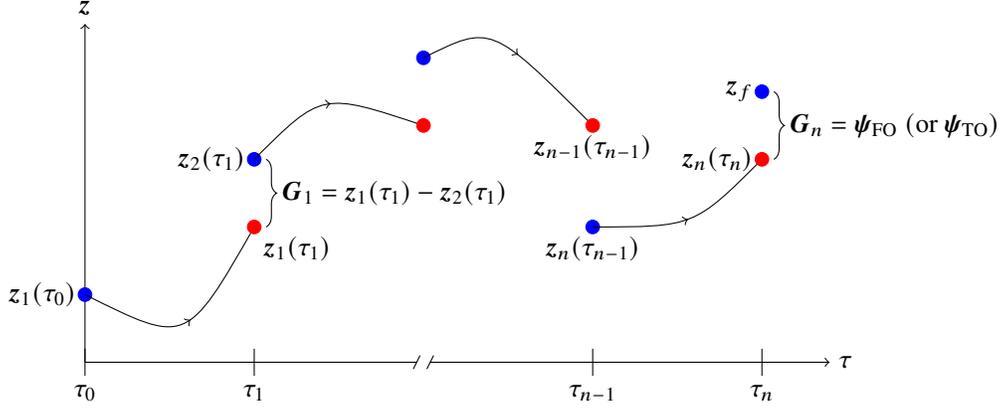

Due to the extreme sensitivity of the problem with respect to the small number of unknown variables, we use an indirect multiple-shooting method similar to the one we used in Ref. \cite{nurre_duty-cycle-aware_2023}. Specifically, this sensitivity arises due to the presence of the singularities of the secant penalty functions on active constraint arcs. The indirect multiple shooting is instrumental in solving the 6DOF PDG problem with SOPICs and it is explained in its entirety here for completeness. As depicted in Fig.~\ref{fig: multiple shooting schematic}, the scaled time-horizon of the problem, \(\tau \in [0,1]\), is divided into \(n\) equal-length intervals such that \(\tau_0=0<\tau_1<\dots<\tau_{n-1}<\tau_n=1\). The states and costates on each \(i\)-th segment, i.e., \(\bm{z}_i(\tau)\), must be continuous between each segment. This constraint is defined as, 
\begin{align}
    \bm{G}_i = \bm{z}_i(\tau_i) - \bm{z}_{i+1}(\tau_i), \quad \forall ~ i = 1,\dots,n-1,
\end{align}
where \(\tau_i\) denotes both the end of the \(i\)-th segment and the beginning of the \((i+1)\)-th segment. The additional unknown variables of the problem, when \(n>1\) is considered, become the states and costates at the beginning of each interval except for the first one, i.e., \(\bm{z}_{i,0}\) such that \(\bm{z}_i(\tau_{i-1})=\bm{z}_{i,0}\) for \(i=2,\cdots, n\). This means the number of unknowns and constraint equations is \(15+28(n-1)\).

We collect all the unknown states and costates into a single vector defined as, $\bm{u} = [\bm{u}_1, \bm{u}_2, \cdots, \bm{u}_n ]$, where $\bm{u}_1 = [\bm{\lambda}^\top_{\bm{r},1,0},\bm{\lambda}^\top_{\bm{v},1,0}, \bm{q}^\top_{1,0}, \bm{\lambda}^\top_{\bm{\omega},1,0}, \lambda_{m,1,0} ]^\top$ and $\bm{u}_i = \bm{z}_{i,0}~\forall~i=2,\cdots, n$. The \(n\)-th constraint vector is defined as the constraint vector at final time, i.e., $\bm{G}_n = \bm{\psi}_\text{FO}$ for the minimum-fuel problem and $\bm{G}_n = \bm{\psi}_\text{TO}$ for the minimum-time problem. We collect all the constraint vectors into the shooting function defined as, $\bm{G}(t_f,\bm{u};\bm{\rho}) = \begin{bmatrix} \bm{G}_i ~ \forall ~ i=1, \cdots, n \end{bmatrix} = \bm{0}$, where the unknown parameter \(t_f\) is additionally an input to the shooting function. The Jacobian of the shooting function with respect to the unknowns has to be calculated. Each row of the Jacobian submatrix with respect to \(t_f\) is defined as,
\begin{equation}\label{eq: jacobian with parameter sensitivity} 
    \frac{\partial \bm{G}_i}{\partial t_f} = \frac{\partial \bm{G}_i}{\partial \bm{z}_i(\tau_i)}\frac{\partial \bm{z}_i(\tau_i)}{\partial t_f}, \quad ~ \forall ~ i = 1, \cdots, n,
\end{equation}
and each row of the Jacobian submatrix with respect to \(\bm{u}\) (i.e., the unknown states and costates at the beginning of each interval) is defined as,
\begin{equation} \label{eq: jacobian with state sensitivity}
    \frac{\partial \bm{G}_i}{\partial \bm{u}_i} = \frac{\partial \bm{G}_i}{\partial \bm{z}_i(\tau_i)}\frac{\partial \bm{z}_i(\tau_i)}{\partial \bm{z}_i(\tau_{i-1})}\frac{\partial \bm{z}_i(\tau_{i-1})}{\partial \bm{u}_i}, ~\quad \forall ~ i = 1, \cdots, n,
\end{equation}
and
\begin{equation}
    \frac{\partial \bm{G}_{i-1}}{\partial \bm{u}_i} = \frac{\partial \bm{G}_{i-1}}{\partial \bm{z}_i(\tau_{i-1})}\frac{\partial \bm{z}_i(\tau_{i-1})}{\partial \bm{u}_i}, ~\quad \forall ~ i = 2 \dots n.
\end{equation}

To show the sparsity pattern of the Jacobian, the \((15 + 28(n-1))\times(15 + 28(n-1))\) Jacobian can be written as,
\begin{equation} \label{eq: bvp jacobian}
    \begin{bmatrix}
        \frac{\partial \bm{G}}{\partial t_f} & \frac{\partial \bm{G}}{\partial \bm{u}}
    \end{bmatrix} = 
    \begin{bmatrix}
        \frac{\partial \bm{G}_1}{\partial t_f} & \frac{\partial \bm{G}_1}{\partial \bm{u}_1} & \frac{\partial \bm{G}_1}{\partial \bm{u}_2} & \bm{0} & \bm{0} & \cdots & \bm{0} \\
        \frac{\partial \bm{G}_2}{\partial t_f} & \bm{0} & \frac{\partial \bm{G}_2}{\partial \bm{u}_2} & \frac{\partial \bm{G}_2}{\partial \bm{u}_3} & \bm{0} & \cdots & \bm{0} \\
        \frac{\partial \bm{G}_3}{\partial t_f} & \bm{0} & \bm{0} & \frac{\partial \bm{G}_3}{\partial \bm{u}_3} & \frac{\partial \bm{G}_3}{\partial \bm{u}_4} & \cdots & \bm{0} \\
        \vdots & \vdots & \vdots & \vdots & \ddots & \ddots & \vdots \\
        \frac{\partial \bm{G}_{n-1}}{\partial t_f} & \bm{0} & \bm{0} & \bm{0} & \cdots & \frac{\partial \bm{G}_{n-1}}{\partial \bm{u}_{n-1}} & \frac{\partial \bm{G}_{n-1}}{\partial \bm{u}_n} \\
        \frac{\partial \bm{G}_n}{\partial t_f} & \bm{0} & \bm{0} & \bm{0} & \cdots & \bm{0} & \frac{\partial \bm{G}_n}{\partial \bm{u}_n}
    \end{bmatrix}_{(15 + 28(n-1))\times(15 + 28(n-1)) }.
\end{equation}

To determine the states and costates at the end of each interval, i.e., \(\bm{z}_i(\tau_i)\) for \(i=1,\cdots, n\), the following IVP is solved 
\begin{align} \label{eq: shooting problem ivp}
    \dot{\bm{z}}_i & = \bm{f}\left(\tau,\bm{z};t_f,\bm{\rho}\right), & \bm{z}_i(\tau_{i-1}) & = \bm{z}_{i,0}, & \tau\in[\tau_{i-1},\tau_i].
\end{align}

The sensitivity matrices in Eqs.~\eqref{eq: jacobian with parameter sensitivity} and \eqref{eq: jacobian with state sensitivity} can be defined as,
\begin{align}
    \bm{\Phi}_i(\tau_i,\tau_{i-1}) & = \frac{\partial \bm{z}_i(\tau_i)}{\partial \bm{z}_i(\tau_{i-1})}, &
    \bm{\Psi}_i(\tau_i,\tau_{i-1}) & = \frac{\partial \bm{z}_i(\tau_i)}{\partial t_f},
\end{align}
and are calculated by propagating their differential equations,
\begin{align}
    \dot{\bm{\Phi}}_i(\tau_i,\tau_{i-1}) & = \frac{\partial \bm{f}\left(\tau,\bm{z};t_f,\bm{\rho}\right)}{\partial \bm{z}_i(\tau)}\bm{\Phi}_i(\tau,\tau_{i-1}), & \bm{\Phi}(\tau_{i-1},\tau_{i-1}) & = \bm{I}_{14\times 14}, \\
    \dot{\bm{\Psi}}_i(\tau_i,\tau_{i-1}) & = \frac{\partial \bm{f}\left(\tau,\bm{z};t_f,\bm{\rho}\right)}{\partial \bm{z}_i(\tau)}\bm{\Psi}_i(\tau,\tau_{i-1}) + \frac{\partial \bm{f}\left(\tau,\bm{z};t_f,\bm{\rho}\right)}{\partial t_f}, & \bm{\Psi}(\tau_{i-1},\tau_{i-1}) & = \bm{0}_{14\times 1},
\end{align}
simultaneously along with Eq.~\eqref{eq: shooting problem ivp}.

The Jacobian function, \(\frac{\partial \bm{f}\left(\tau,\bm{z};t_f,\bm{\rho}\right)}{\partial \bm{z}_i(\tau)}\), is calculated using AD of CasADi. The Jacobian \(\frac{\partial \bm{f}\left(\tau,\bm{z};t_f,\bm{\rho}\right)}{\partial t_f}\) is trivially calculated (see Eq.~\eqref{eq: bvp system}). The nonlinear shooting problems are solved using MATLAB's \verb|fsolve|. Because the BVPs demonstrates a large condition number in their Jacobian, the Levenberg-Marquardt algorithm with the option \verb|ScaleProblem| set to \verb|jacobian| was found to be the most effective algorithm. A \verb|FunctionTolerance| and \verb|StepTolerance| of \(1.0 \times 10^{-10}\) and \(1.0 \times 10^{-14}\) were used, respectively. The IVPs are integrated using MATLAB's \verb|ode113| with absolute and relative integration tolerances set to \(1.0 \times 10^{-12}\).

\section{Results} \label{sec: results}
For the 6DOG PDG problems, the parameters and boundary conditions in Tables \ref{tab: problem parameters} and \ref{tab: boundary conditions} are taken from Ref. \cite{sagliano_six-degree--freedom_2024}. Non-dimensionalized units for length, time, and mass are used which are denoted as LU, TU, and MU, respectively. A force unit is also defined as FU = MU\(\cdot\)LU/TU\(^2\). We emphasize that minor changes are applied to final position values to avoid singularities due to \(r_x=r_y=r_z=0\) in $S_2$ in Eq.~\eqref{eq: new control constraints}. The final vertical position, \(r_z(t_f)\), is set equal to the gimbal moment arm length, \(\|\bm{r}_T\|\), which makes more practical sense than having the rocket's center of mass at ground level and, thus, its gimbal point under the surface of the Earth. Only the final vertical position needs to be offset because we implement this solution methodology in MATLAB, which calculates \(\tan^{-1}{\left(1/0\right)} = \frac{\pi}{2}\) (to machine precision) and \(\tan^{-1}{\left(0/0\right)}=\verb|NaN|\). 
\begin{remark}
    This is likely the reason for the spike at the end of the glideslope angle profile in Figure 8 of Ref. \cite{sagliano_six-degree--freedom_2024}. Our indirect method can converge when the final position boundary condition is \((0,0,0)\) but it also exhibits the same spike at the end of the glideslope angle profile.
\end{remark} 
We also solve the problems using the student version of DIDO 7.5.6 \cite{ross_enhancements_2020} for validation of our solutions. We emphasize that the
solutions are obtained independent of each other (i.e., DIDO’s solution is not used to initialize the indirect method).

\begin{table}[!tb]
    \centering
    \caption{Problem parameters (Ref. \cite{sagliano_six-degree--freedom_2024}).}
    \label{tab: problem parameters}
    \begin{tabular}{lc}
        \hline
        Parameter & Value [unit] \\ \hline \hline
        Wet mass, \(m_\text{wet}\) & 2 [MU] \\ \hline
        Dry mass, \(m_\text{dry}\) & 1 [MU] \\ \hline
        Drag coefficient, \(C_D\) & 0.1 \\ \hline
        Atmospheric density, \(\rho\) & 1 [MU/LU\(^2\)] \\ \hline
        Reference surface, \(S\) & 0.5 [LU\(^2\)] \\ \hline
        Minimum thrust, \(T_\text{min}\) & 1 [FU]\\ \hline
        Maximum thrust, \(T_\text{max}\) & 5 [FU] \\ \hline
        Specific impulse, \(I_\text{sp}\) & 294.2 [TU] \\ \hline
        Sea-level gravity, \(g_0\) & 1 [LU/TU\(^2\)] \\ \hline
        Inertia matrix, \(\bm{J}\) & \(0.01\cdot \bm{I}_{3\times3}\) [MU\(\cdot\)LU\(^2\)] \\ \hline
        Gimbal point location, \(\bm{r}_T\) & \(\left[0,0,-0.01\right]^\top\) [LU] \\ \hline
        Minimum glideslope angle, \(\gamma_\text{min}\) & 20 [deg] \\ \hline
        Maximum tilt angle, \(\theta_\text{max}\) & 90 [deg] \\ \hline
        Maximum angular rate, \(\omega_\text{max}\) & 60 [deg/TU] \\ \hline
        Maximum gimbal angle, \(\delta_\text{max}\) & 20 [deg] \\ \hline
    \end{tabular}
\end{table}

\begin{table}[!tb]
    \centering
    \caption{Boundary conditions (Ref. \cite{sagliano_six-degree--freedom_2024}) with slight modifications to avoid certain singularities.}
    \label{tab: boundary conditions}
    \begin{tabular}{lc}
        \hline
        Boundary Condition & Value [unit] \\ \hline \hline
        Initial position, \(\bm{r}_0\) & \(\left[0.5,4,4\right]^\top\) [LU] \\ \hline
        Initial velocity, \(\bm{v}_0\) & \(\left[0,-4,0\right]^\top\) [LU/TU] \\ \hline
        Initial angular velocity, \(\bm{\omega}_0\) & \(\left[0,0,0\right]^\top\) [1/TU] \\ \hline
        Initial mass, \(m_0\) & \(m_\text{wet}\) [MU] \\ \hline
        Final position, \(\bm{r}_f\) & \(\left[0,0,0.01\right]^\top\) [LU] \\ \hline
        Final velocity, \(\bm{v}_f\) & \(\left[0,0,0\right]^\top\) [LU/TU] \\ \hline
        Final orientation, \(\bm{q}_f\) & \(\left[0,0,0.01,1\right]^\top\) [-] \\ \hline
        Final angular velocity, \(\bm{\omega}_f\) & \(\left[0,0,0\right]^\top\) [1/TU] \\ \hline
    \end{tabular}
\end{table}

\subsection{Fuel-Optimal Solution}
To solve the fuel-optimal problem through continuation, the problem was first solved with the smoothing parameters set to \(\rho_T=0.01\), \(\rho_\delta=0.1\), \(\rho_\omega=0\), \(\rho_\gamma=0\), and \(\rho_\theta=0\), which represents the problem with no SOPICs. This could be solved with a single-shooting method, i.e., \(n=1\) segments. This solution had an initial tilt angle that violated \(\theta_\text{max}\) so \(\theta_\text{max}=95\) deg was set so that the solution was interior to the constraint. Continuation was then performed over all the parameters (except \(\rho_\omega=0\) and \(\rho_\gamma=0\) since the maximum angular velocity and minimum glideslope angle constraints are never active) and \(\theta_\text{max}\) until their values were \(\rho_T=10^{-7}\), \(\rho_\delta=10^{-4}\), \(\rho_\theta=10^{-12}\), and \(\theta_\text{max}=90\) deg. 
\begin{figure}[!htbp]
    \centering
    \includegraphics[width=0.7\textwidth]{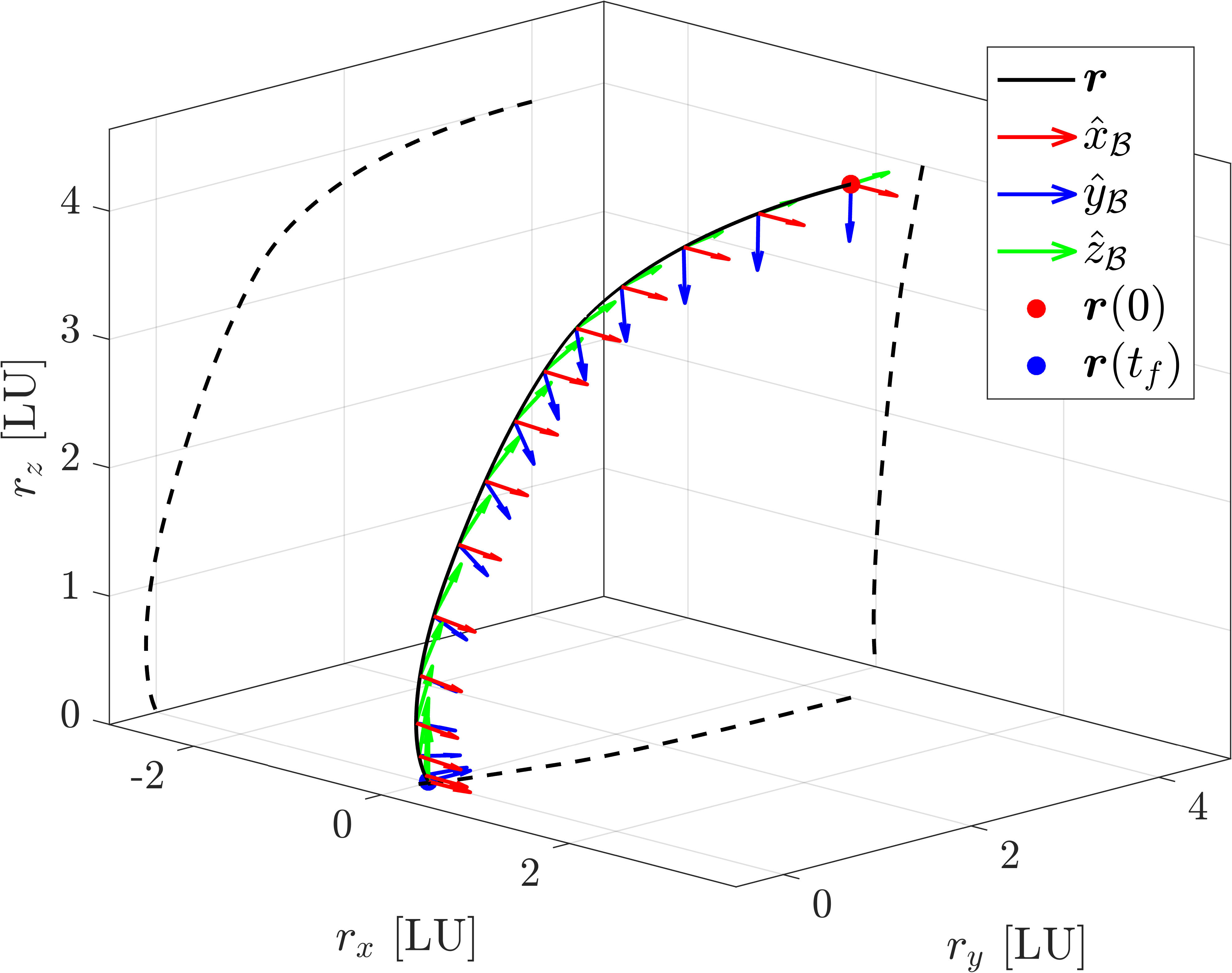}
    \caption{Fuel-optimal: trajectory solution and its projection on different inertial planes.}
    \label{fig: p1 trajectory}
\end{figure}
The number of shooting segments had to be increased to \(n=5\) to solve this problem. The final mass was \(m(t_f)=1.95382\) MU and the time of flight found was \(t_f=3.72457\) TU. The final mass and time of flight from the DIDO solution were \(m(t_f)=1.95381\) MU and \(t_f=3.72677\) TU, respectively. The final mass for the solution obtained in Ref. \cite{sagliano_six-degree--freedom_2024} was \(m(t_f)=1.9523\) MU \footnote{Through personal communications with Dr. Marco Sagliano.}, which is based on a pseudo-spectral convex optimization method.

\begin{figure}[]
    \centering
    \begin{subfigure}[]{0.49\textwidth}
        \centering
        \includegraphics[width=\textwidth]{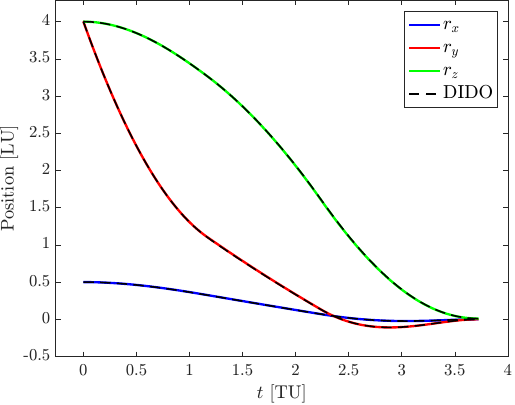}
        \caption{Position}
        \label{fig: p1 states position}
    \end{subfigure}%
    ~ 
    \begin{subfigure}[]{0.49\textwidth}
       \centering
       \includegraphics[width=\linewidth]{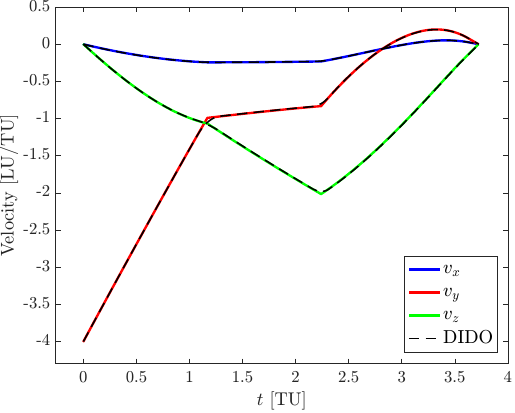}
       \caption{Velocity}
       \label{fig: p1 states velocity}
    \end{subfigure}%
    \\ \vspace{1mm}
    \begin{subfigure}[]{0.49\textwidth}
        \centering
        \includegraphics[width=\textwidth]{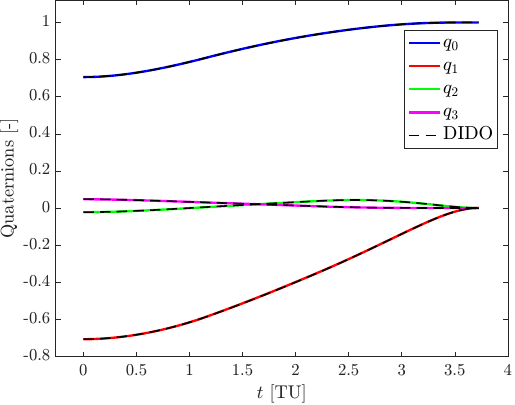}
        \caption{Quaternions}
        \label{fig: p1 states quaternions}
    \end{subfigure}%
    ~
    \begin{subfigure}[]{0.49\textwidth}
        \centering
        \includegraphics[width=\textwidth]{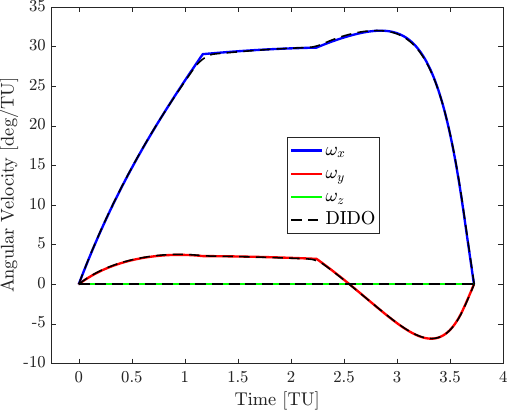}
        \caption{Angular Velocity}
        \label{fig: p1 states angular velocity}
    \end{subfigure}
    \caption{Fuel-optimal: state time histories.}
    \label{fig: p1 states}
\end{figure}

Figure \ref{fig: p1 trajectory} shows the fuel-optimal trajectory and its projections on different planes, with the body axes overlaid to show the rocket's orientation at various points along it. Figure \ref{fig: p1 states} shows the time histories of the states obtained using the indirect method and DIDO. The components of the quaternions are plotted in Figure \ref{fig: p1 states}.

\begin{figure}[]
    \centering
    \includegraphics[width=0.9\textwidth]{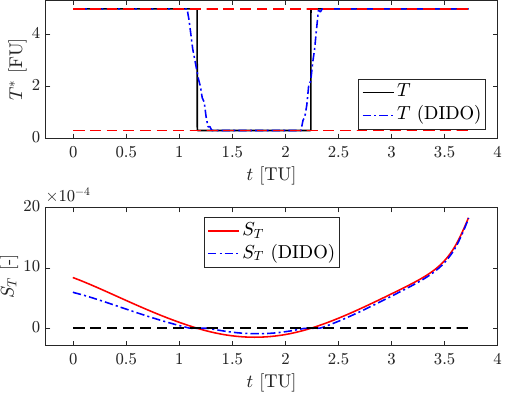}
    \caption{Fuel-optimal: thrust control and switching function profiles.}
    \label{fig: p1 thrust}
\end{figure}

\begin{figure}[]
    \centering
    \begin{subfigure}[]{0.49\textwidth}
        \centering
        \includegraphics[width=\textwidth]{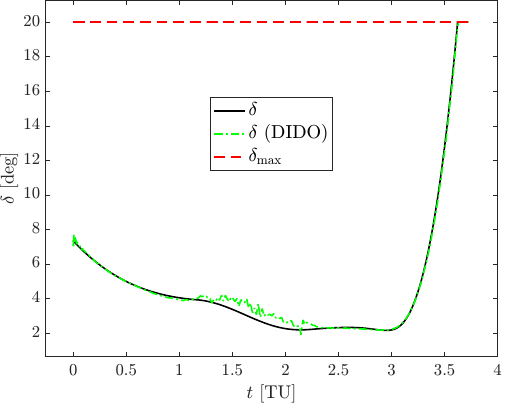}
        \caption{Gimbal angle profile.}
        \label{fig: gimbal angle}
    \end{subfigure}%
    ~
    \begin{subfigure}[]{0.49\textwidth}
        \centering
        \includegraphics[width=\textwidth]{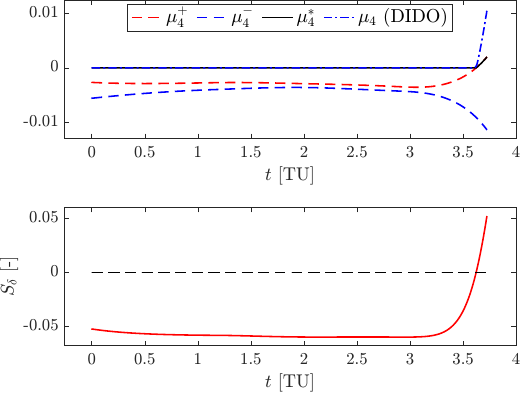}
        \caption{Lagrange multipliers and switching function profiles.}
        \label{fig: Lagrange gimbal constraint}
    \end{subfigure}
    \caption{Fuel-optimal: gimbal angle, and its inequality Lagrange multipliers and switching function profiles.}
    \label{fig: fo steering control}
\end{figure}

Figures \ref{fig: p1 thrust} shows the thrust magnitude control profile along with its switching function. The thrust switching function from DIDO is the Lagrange multiplier associated with this inequality constraint. Figure \ref{fig: gimbal angle} shows that the gimbal angle constraint is saturated at the end of the maneuver. Figure \ref{fig: Lagrange gimbal constraint} shows the behavior of \(\tilde{\mu}_4^*\) and \(S_\delta\) accompanying this constraint activation. The indirect solution agrees well with DIDO and is similar to the solution presented in Ref.~\cite{sagliano_six-degree--freedom_2024} (other than discrepancies in the middle of the time horizon in the gimbal angle in Ref.~\cite{sagliano_six-degree--freedom_2024}). 

\begin{figure}[]
    \centering
    \begin{subfigure}[]{0.49\textwidth}
        \centering
        \includegraphics[width=\textwidth]{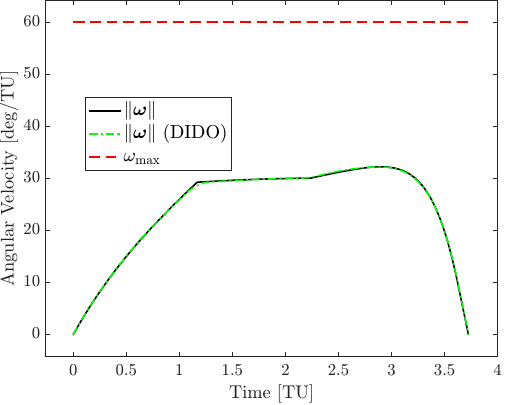}
        \caption{Angular velocity magnitude vs time.}
        \label{fig: p1 angular velocity magnitude}
    \end{subfigure}%
    ~ 
    \begin{subfigure}[]{0.49\textwidth}
        \centering
        \includegraphics[width=\textwidth]{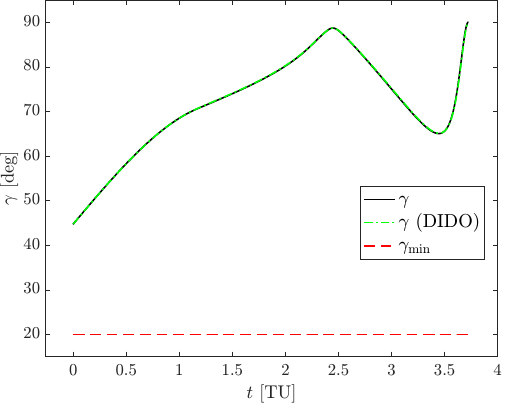}
        \caption{Glideslope angle vs. time.}
        \label{fig: p1 glideslope angle}
    \end{subfigure}
    \caption{Fuel-optimal: angular velocity magnitude and glideslope angle profiles.}
    \label{fig: fo angular velocity and glideslope}
\end{figure}

\begin{figure}[]
    \centering
    \begin{subfigure}[]{0.49\textwidth}
        \centering
        \includegraphics[width=\textwidth]{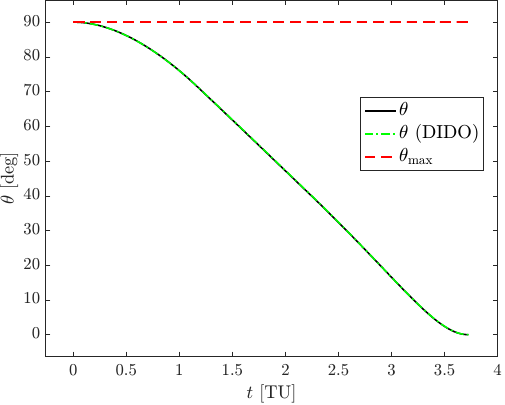}
        \caption{Tilt angle vs. time.}
        \label{fig: p1 tilt angle}
    \end{subfigure}%
    ~
    \begin{subfigure}[]{0.49\textwidth}
        \centering
        \includegraphics[width=\textwidth]{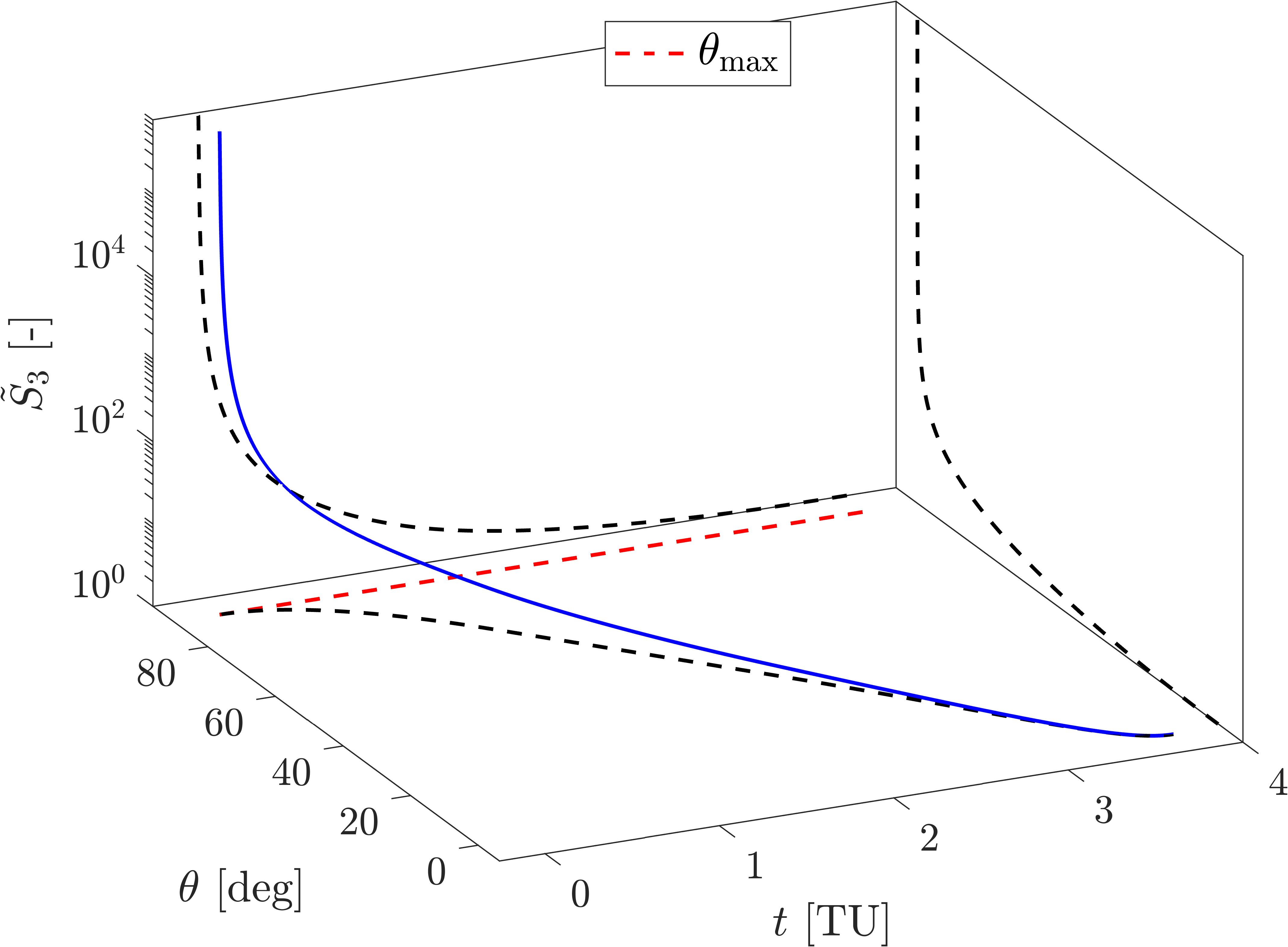}
        \caption{Tilt angle constraint secant penalty function profile.}
        \label{fig: p1 tilt angle constraint}
    \end{subfigure}
    \caption{Fuel-optimal: tilt angle constraint.}
    \label{fig: fo tilt constraints}
\end{figure}

Inspecting Figures \ref{fig: p1 angular velocity magnitude}, \ref{fig: p1 glideslope angle}, and \ref{fig: p1 tilt angle}, it can be seen the tilt angle constraint was the only SOPIC that became active and only for an instant at the very beginning of the solution. The initial optimal orientation was found to be \(\bm{q}(0)=\left[0.705530156561408, -0.706747429691534, -0.0225333659710576, 0.0471965705287370\right]^\top\), which results in an initial tilt angle of \(\theta(0)=89.99996\) deg. The trend was observed that as \(\rho_\theta \rightarrow 0\) then \(\theta(0)\rightarrow 90\) deg. Figure \ref{fig: p1 tilt angle constraint} shows the time history of the secant penalty function, \(\tilde{S}_{3}\) (Eq.~\eqref{eq: secant penalty functions}), as a function of the tilt angle.

\begin{figure}[]
    \centering
    \begin{subfigure}[]{0.49\textwidth}
        \centering
        \includegraphics[width=\textwidth]{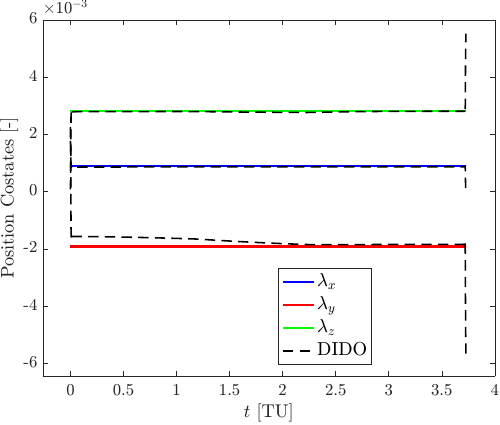}
        \caption{Position Costates}
        \label{fig: p1 costates position}
    \end{subfigure}%
    ~ 
    \begin{subfigure}[]{0.49\textwidth}
        \centering
        \includegraphics[width=\textwidth]{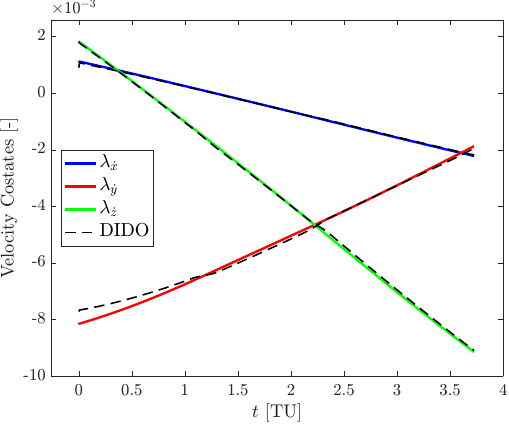}
        \caption{Velocity Costates}
        \label{fig: p1 costates velocity}
    \end{subfigure}%
    \\ \vspace{1mm}
    \begin{subfigure}[]{0.49\textwidth}
        \centering
        \includegraphics[width=\textwidth]{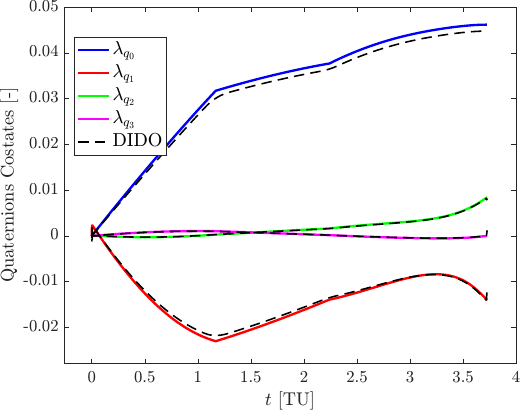}
        \caption{Quaternion Costates}
        \label{fig: p1 costates quaternions}
    \end{subfigure}%
    ~
    \begin{subfigure}[]{0.49\textwidth}
        \centering
        \includegraphics[width=\textwidth]{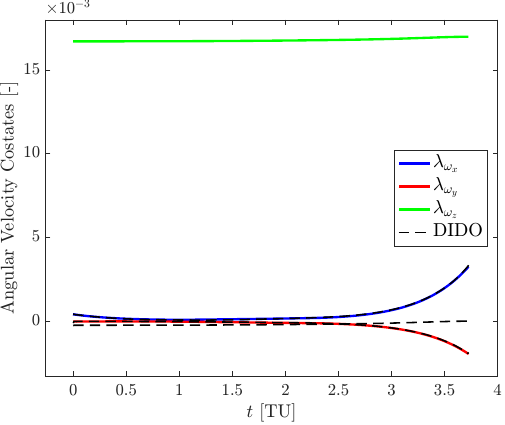}
        \caption{Angular Velocity Costates}
        \label{fig: p1 costates angular velocity}
    \end{subfigure}
    \caption{Fuel-optimal: costate time histories.}
    \label{fig: p1 costates}
\end{figure}

\begin{figure}[]
    \centering
    \begin{subfigure}[]{0.49\textwidth}
        \centering
        \includegraphics[width=\textwidth]{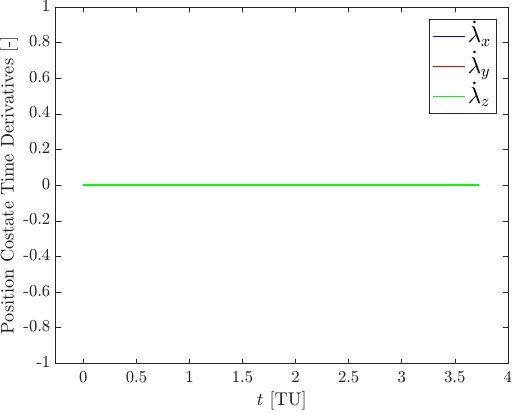}
        \caption{Position Costate Time Derivatives}
        \label{fig: p1 costates time derivatives position}
    \end{subfigure}%
    ~ 
    \begin{subfigure}[]{0.49\textwidth}
        \centering
        \includegraphics[width=\textwidth]{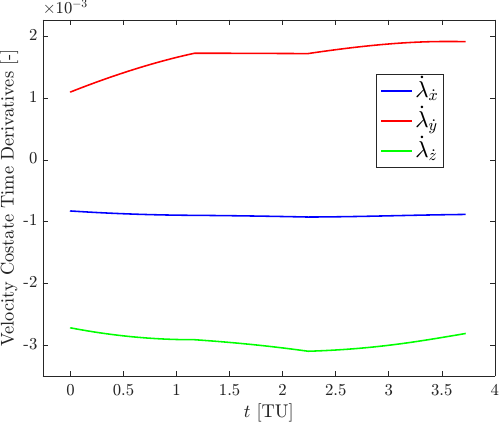}
        \caption{Velocity Costate Time Derivatives}
        \label{fig: p1 costates time derivatives velocity}
    \end{subfigure}%
    \\ \vspace{1mm}
    \begin{subfigure}[]{0.49\textwidth}
        \centering
        \includegraphics[width=\textwidth]{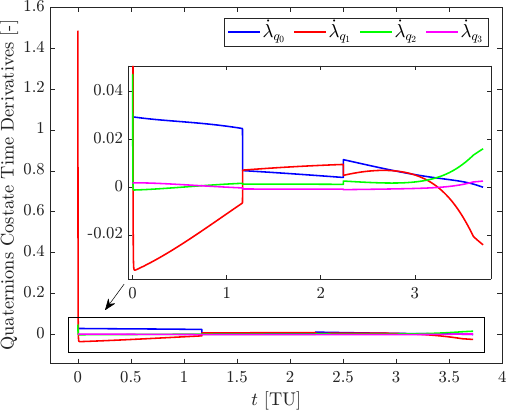}
        \caption{Quaternion Costate Time Derivatives}
        \label{fig: p1 costates time derivatives quaternions}
    \end{subfigure}%
    ~
    \begin{subfigure}[]{0.49\textwidth}
        \centering
        \includegraphics[width=\textwidth]{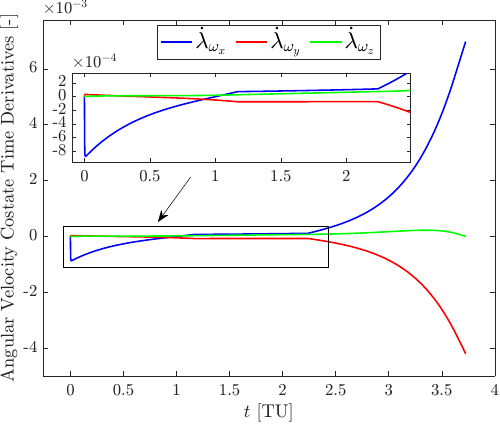}
        \caption{Angular Velocity Costate Time Derivatives}
        \label{fig: p1 costates time derivatives angular velocity}
    \end{subfigure}
    \caption{Fuel-optimal: time histories of the time-derivative of the costates.}
    \label{fig: p1 costates derivatives}
\end{figure}

\begin{figure}[]
    \centering
    \begin{subfigure}[]{0.49\textwidth}
        \centering
        \includegraphics[width=\textwidth]{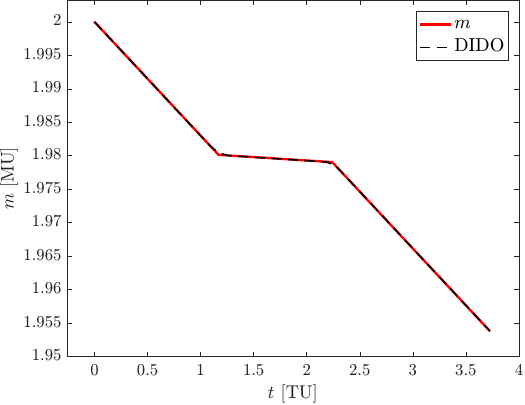}
        \caption{Mass}
        \label{fig: p1 subplot mass}
    \end{subfigure}%
    ~ 
    \begin{subfigure}[]{0.49\textwidth}
        \centering
        \includegraphics[width=\textwidth]{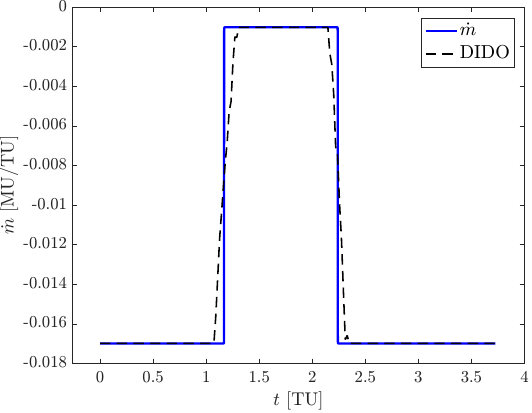}
        \caption{Mass Time Derivative}
        \label{fig: p1 mass time derivative}
    \end{subfigure}%
    \\ \vspace{1mm}
    \begin{subfigure}[]{0.49\textwidth}
        \centering
        \includegraphics[width=\textwidth]{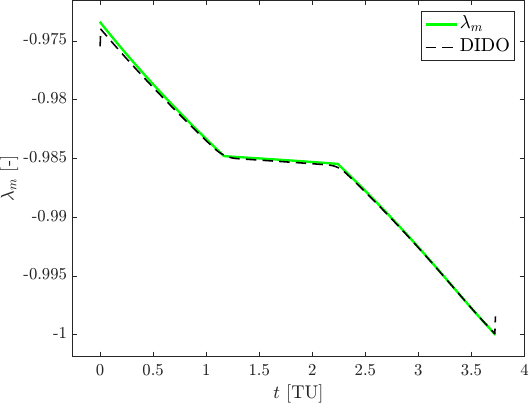}
        \caption{Mass Costate}
        \label{fig: p1 mass costate}
    \end{subfigure}%
    ~
    \begin{subfigure}[]{0.49\textwidth}
        \centering
        \includegraphics[width=\textwidth]{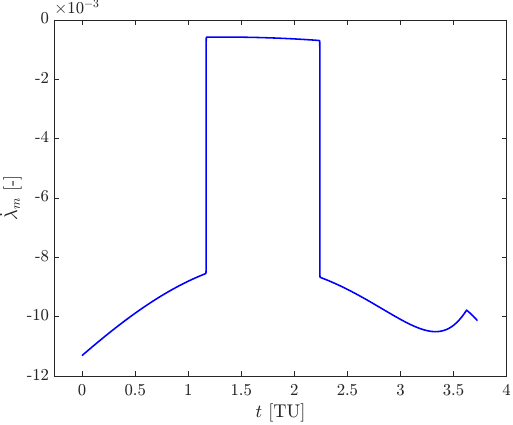}
        \caption{Mass Costate Time Derivative}
        \label{fig: p1 mass costate time derivative}
    \end{subfigure}
    \caption{Fuel-optimal: mass, mass flow, mass costate, and mass costate time-derivative time histories.}
    \label{fig: p1 mass}
\end{figure}

\begin{figure}
    \centering
    \includegraphics[width=0.4\textwidth]{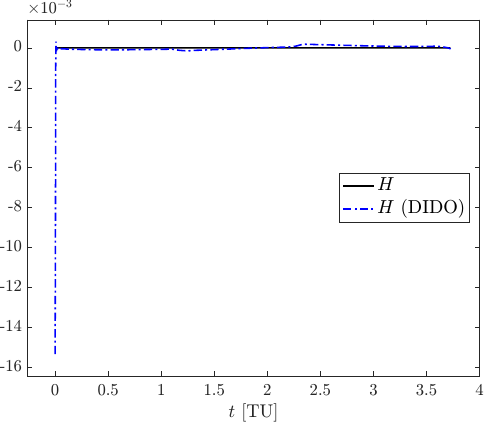}
    \caption{Fuel-optimal: Hamiltonian profile.}
    \label{fig: p1 hamiltonian}
\end{figure}

Figures \ref{fig: p1 costates} and \ref{fig: p1 costates derivatives} show the time histories of the costates and their time-derivatives. Note that DIDO does not directly provide costate time-derivative in its solution. Jump discontinuities can be observed in the indirect method solution at initial time in the quaternion costates and the quaternion and angular velocity costate time-derivatives due to the activation of the second-order tilt angle state-path constraint. There is a discrepancy between the indirect and DIDO solutions in the costate for the angular velocity in the \(\hat{\bm{z}}_\mathcal{B}\) direction, \(\lambda_{\omega_z}\). This is likely due to the roll (and thus the longitudinal angular velocity) being uncontrollable in the problem formulation, leading to nonuniqueness of the associated costates. Figure \ref{fig: p1 mass} shows the time histories of the mass, mass flow, mass costate, and mass costate time-derivative. Figure \ref{fig: p1 hamiltonian} shows the Hamiltonian time history, which appears constant and equal to 0, satisfying the stationary condition. In the DIDO solution there is a spike at the beginning, likely having to do with the tilt angle constraint activation, but remains more or less constant after that. 

\subsection{Time-Optimal Solution}

Similar to the fuel-optimal problem, the time-optimal problem was first solved with a single-shooting method and with the smoothing parameters initially set to \(\rho_T=0.01\), \(\rho_\delta=0.1\), \(\rho_\omega=0\), \(\rho_\gamma=0\), and \(\rho_\theta=0\), which represents the problem with no SOPICs. This initial solution violated the maximum tilt angle and angular velocity magnitude values in Table \ref{tab: problem parameters}. Thus, these values were adjusted accordingly to ensure the solution was interior to the constraint set before performing continuation on the rest of the parameters. Continuation was performed over all the parameters (except \(\rho_\gamma=0\) since the glideslope angle constraint is never active) along with \(\theta_\text{max}\) and \(\omega_\text{max}\) until their values were \(\rho_T=\rho_\theta=\rho_\omega=1.0 \times 10^{-7}\),~\(\rho_\delta=1.0 \times 10^{-4}\),~\(\theta_\text{max}=90\) deg, and \(\omega_\text{max}=60\) deg/TU. The number of shooting segments had to be increased to $n=500$ to perform this continuation. This higher number of segments was needed due to the longer duration that the SOPICs were active. This made a large interval of the solution very close to the singularities in the secant penalty functions, which drastically increases the sensitivity of the problem. Some type of mesh refinement scheme could be employed to more efficiently choose the number of segments and their temporal length, but this was beyond the scope of this work. Despite the large number of segments, the problem is parallelizable and the IVPs were solved in parallel in MATLAB with \verb|parfor|. The final mass was \(m(t_f)=1.94977\) MU and the time of flight found was \(t_f=3.50453\) TU. The final mass and time of flight from the DIDO solution (solved with 240 nodes) were \(m(t_f)=1.9498\) MU and \(t_f=3.50453\) TU, respectively.

\begin{figure} [!htbp]
    \centering
    \includegraphics[width=0.8\textwidth]{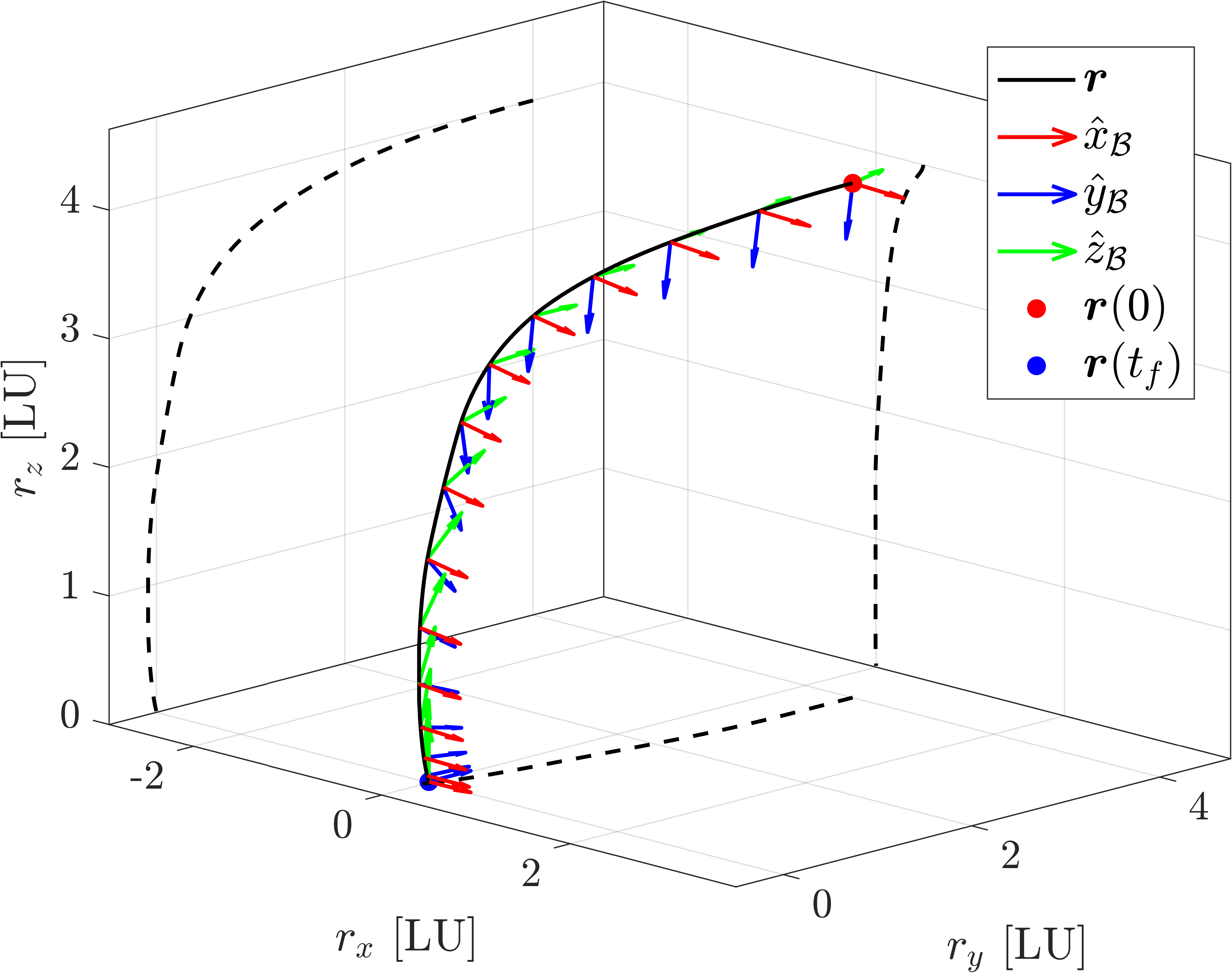}
    \caption{Time-optimal: trajectory and its projection on different inertial planes.}
    \label{fig: p1 mt trajectory}
\end{figure}

\begin{figure}
    \centering
    \resizebox{\linewidth}{!}{%
        \includegraphics{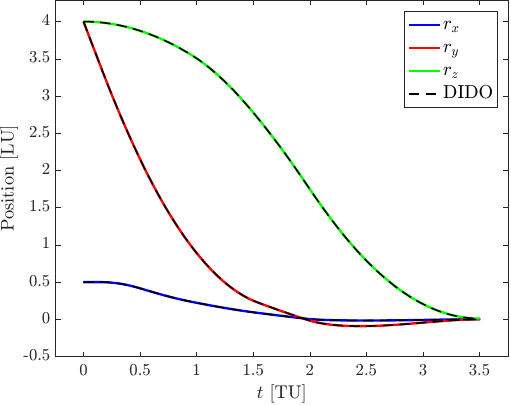}
        \includegraphics{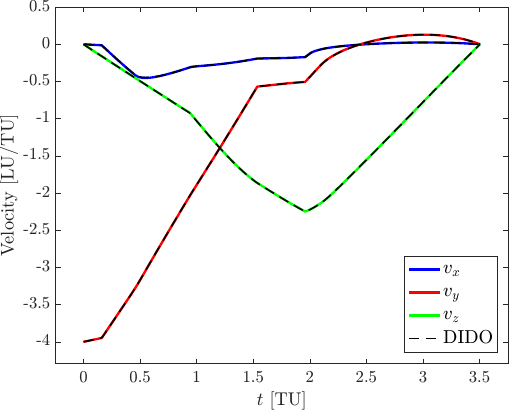}
    }
    \resizebox{\linewidth}{!}{%
        \includegraphics{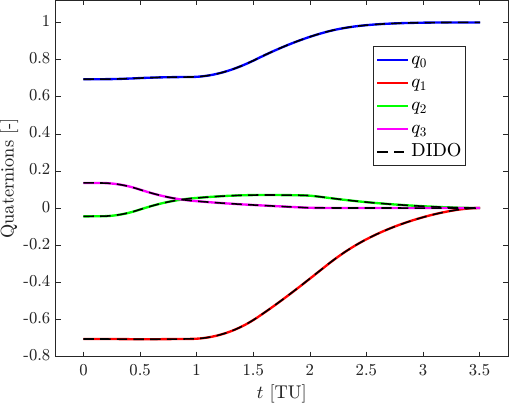}
        \includegraphics{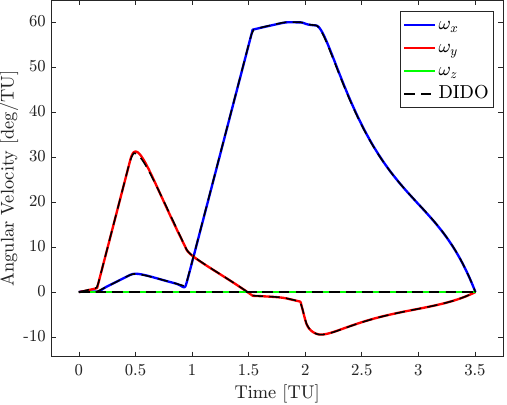}
    }
    \caption{Time-optimal: state time histories.}
    \label{fig: p1 mt states}
\end{figure}

\begin{figure}
    \centering
    \includegraphics[width=0.70\textwidth]{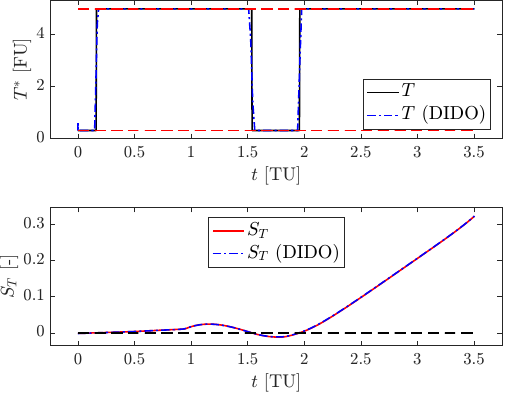}
    \caption{Time-optimal: thrust control and switching function vs. time.}
    \label{fig: p1 mt thrust}
\end{figure}

\begin{figure}
    \centering
    \includegraphics[width=0.75\textwidth]{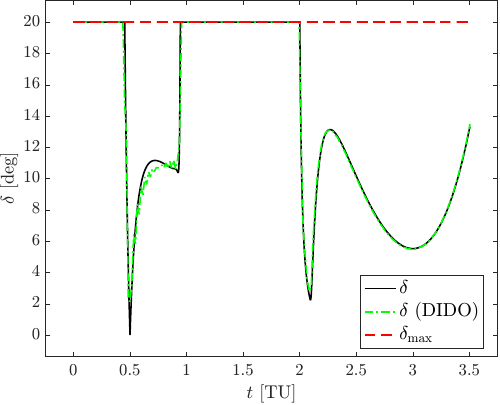}
    \caption{Time-optimal: gimbal angle profile.}
    \label{fig: mt gimbal angle}
\end{figure}

\begin{figure}
    \centering
    \includegraphics[width=0.75\textwidth]{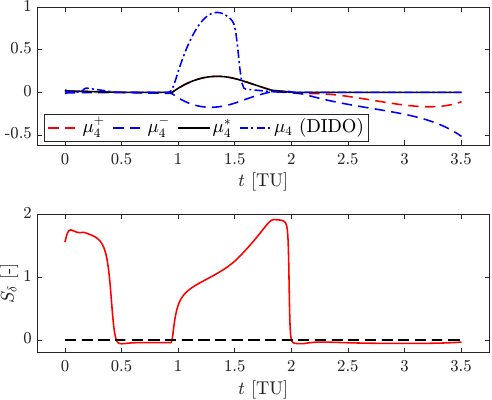}
    \caption{Time-optimal: gimbal angle constraint Lagrange multipliers and switching function profiles.}
    \label{fig: mt Lagrange gimbal constraint}
\end{figure}

\begin{figure}
    \centering
    \includegraphics[width=0.75\textwidth]{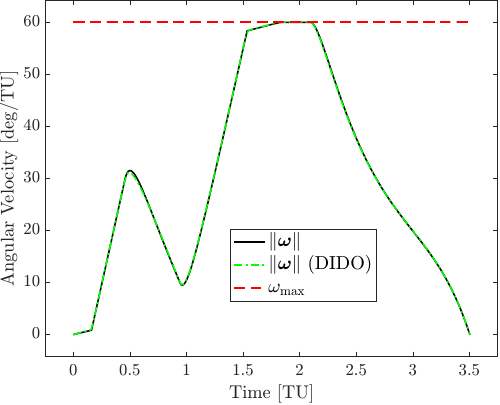}
    \caption{Time-optimal: angular velocity magnitude profile.}
    \label{fig: p1 mt angular velocity magnitude}
\end{figure}

\begin{figure}[]
    \centering
    \includegraphics[width=0.8\textwidth]{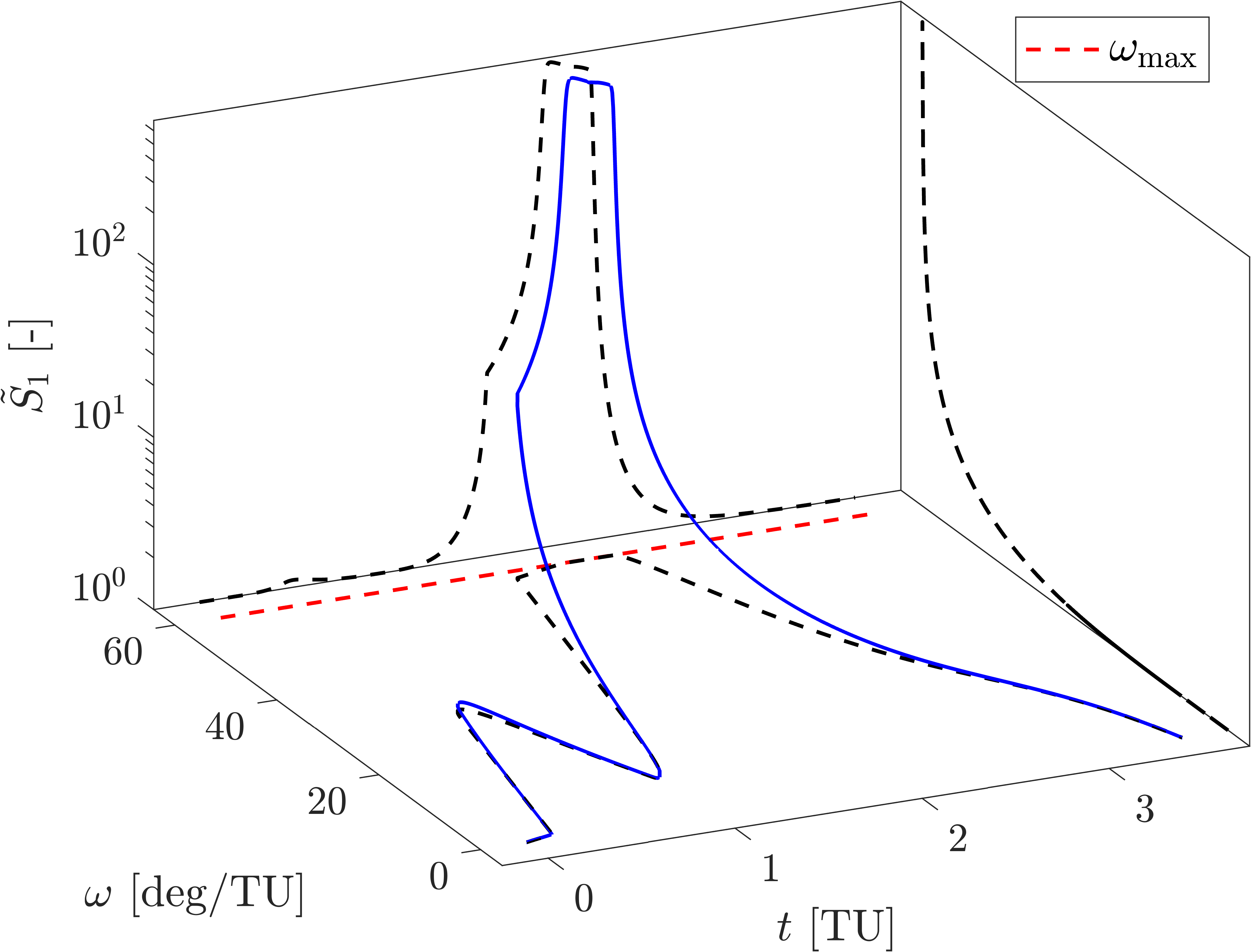}
    \caption{Time-optimal: angular velocity magnitude constraint secant penalty function profile.}
    \label{fig: p1 mt angular velocity magnitude constraint}
\end{figure}

\begin{figure}
    \centering
    \includegraphics[width=0.75\textwidth]{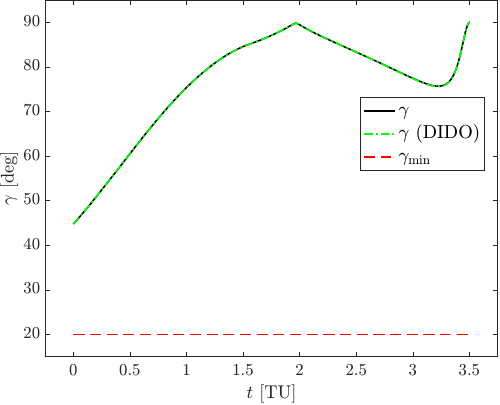}
    \caption{Time-optimal: glideslope angle vs. time.}
    \label{fig: p1 mt glideslope angle}
\end{figure}

\begin{figure}
    \centering
    \includegraphics[width=0.70\textwidth]{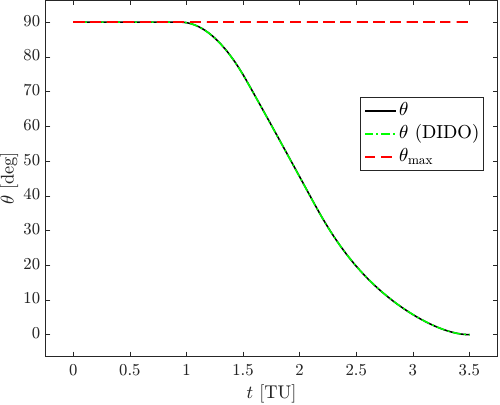}
    \caption{Time-optimal: tilt angle vs. time.}
    \label{fig: p1 mt tilt angle}
\end{figure}

\begin{figure}[]
    \centering
    \includegraphics[width=0.8\textwidth]{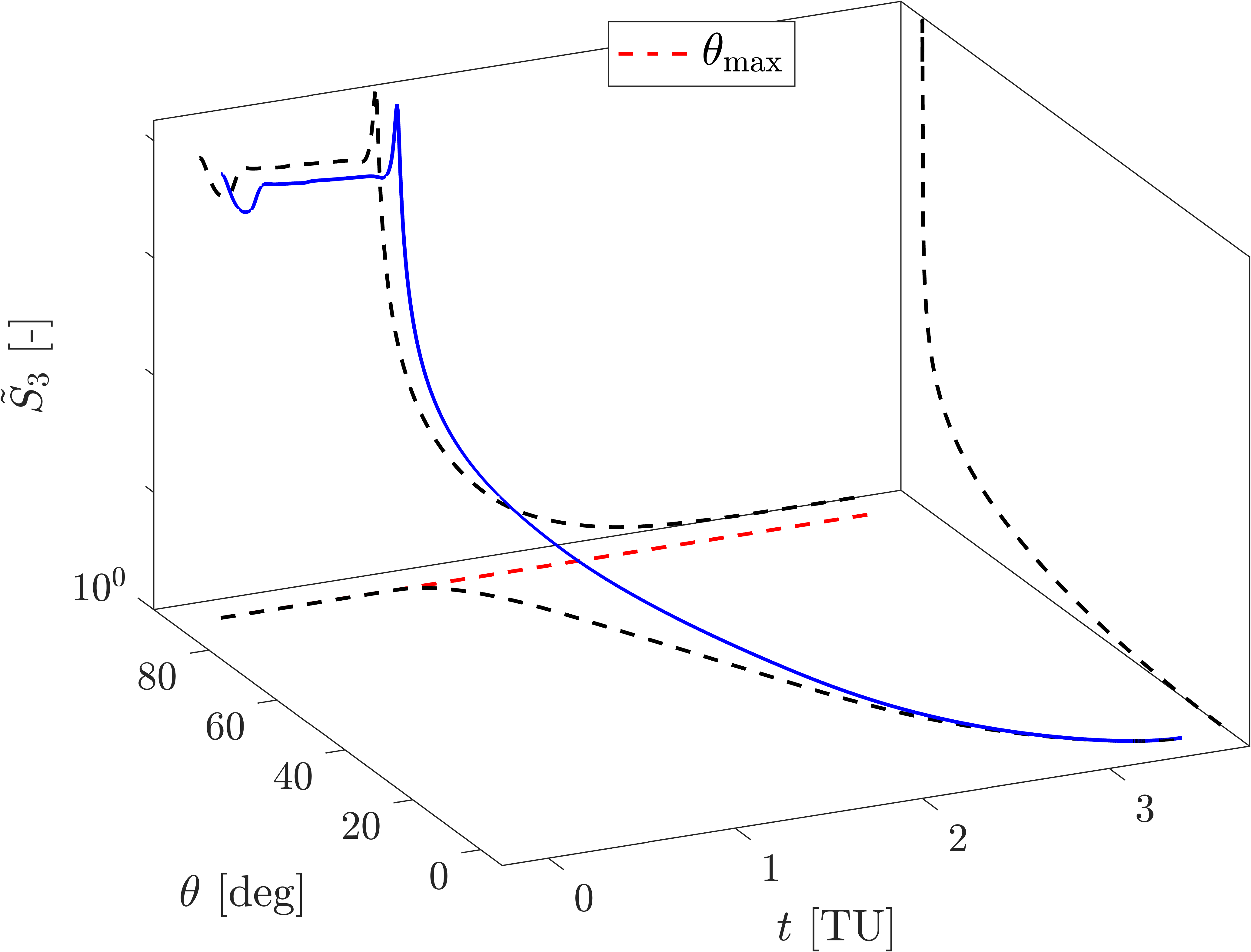}
    \caption{Time-optimal: tilt angle constraint secant penalty function profile.}
    \label{fig: p1 mt tilt angle constraint}
\end{figure}

\begin{figure}
    \centering
    \resizebox{\linewidth}{!}{%
        \includegraphics{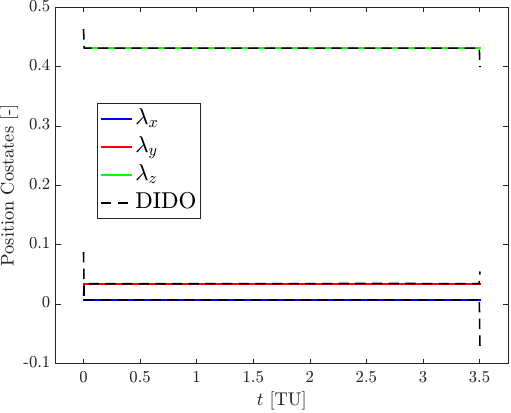}
        \includegraphics{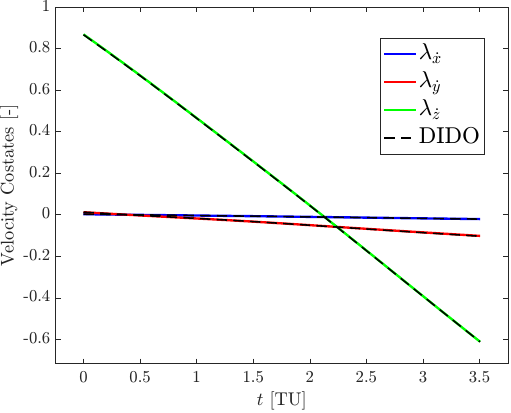}
    }
    \resizebox{\linewidth}{!}{%
        \includegraphics{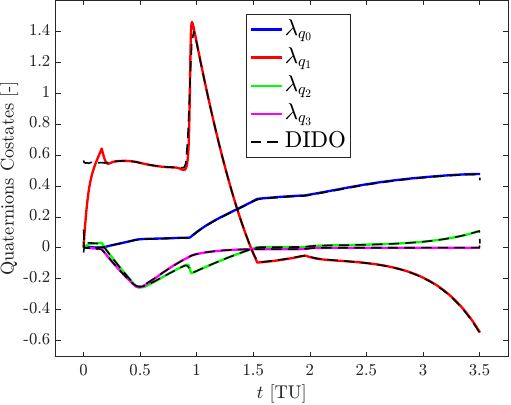}
        \includegraphics{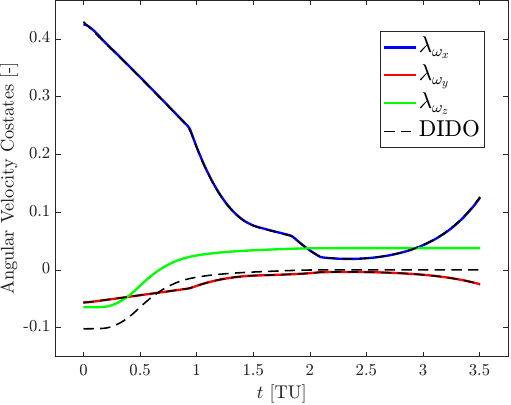}
    }
    \caption{Time-optimal: costate time histories.}
    \label{fig: p1 mt costates}
\end{figure}

\begin{figure}
    \centering
    \resizebox{\linewidth}{!}{%
        \includegraphics{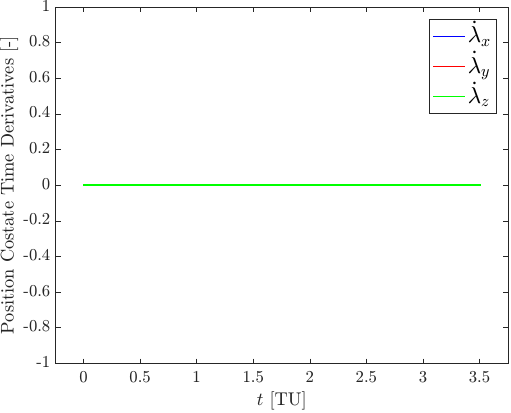}
        \includegraphics{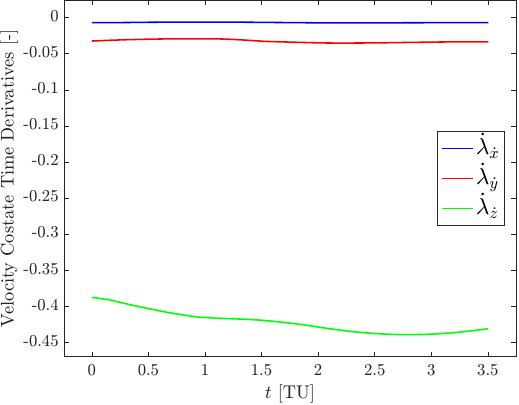}
    }
    \resizebox{\linewidth}{!}{%
        \includegraphics{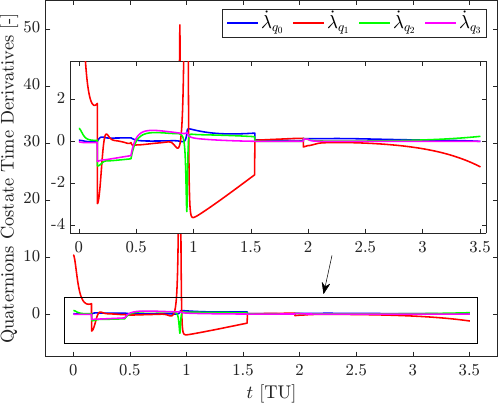}
        \includegraphics{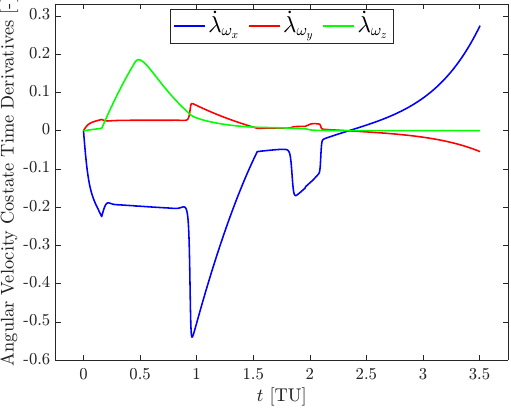}
    }
    \caption{Time-optimal: time histories of the time-derivative of the costates (indirect solution).}
    \label{fig: p1 mt costate derivatives}
\end{figure}

\begin{figure}
    \centering
    \includegraphics[width=0.49\linewidth]{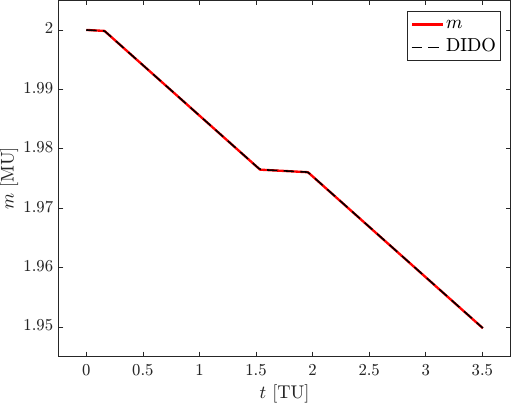}
    \includegraphics[width=0.49\linewidth]{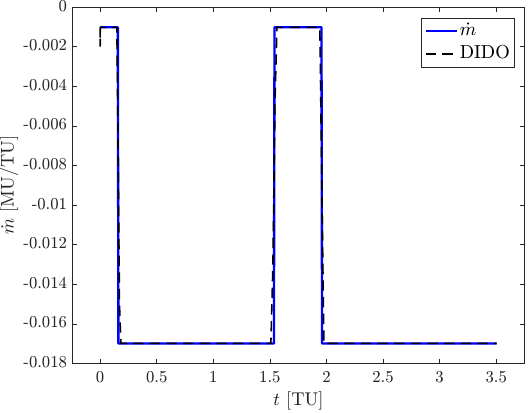} \\
    \includegraphics[width=0.49\linewidth]{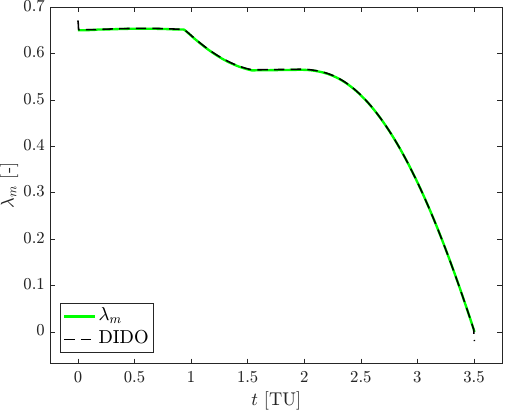}
    \includegraphics[width=0.49\linewidth]{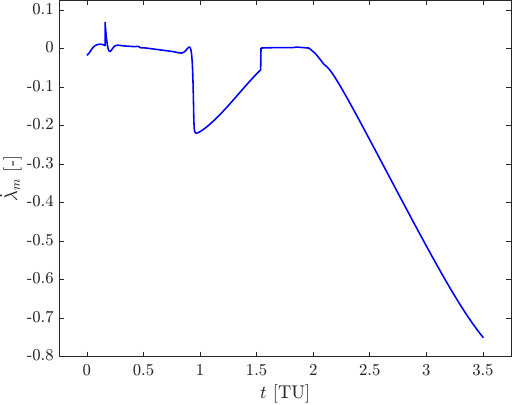}
    \caption{Time-optimal: mass, mass flow, mass costate, and mass costate time-derivative time histories.}
    \label{fig: p1 mt mass}
\end{figure}

\begin{figure}
    \centering
    \includegraphics[width=0.50\textwidth]{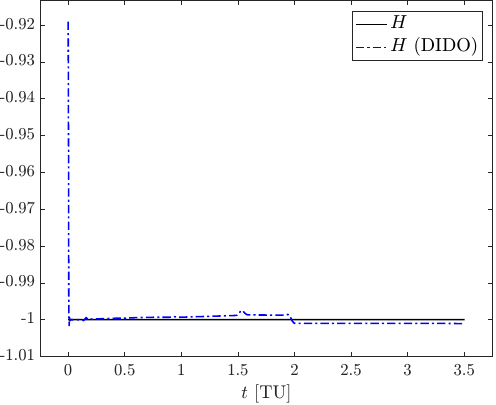}
    \caption{Time-optimal: Hamiltonian profile.}
    \label{fig: mt hamiltonian}
\end{figure}

Figure \ref{fig: p1 mt trajectory} shows the time-optimal trajectory with the body axes overlaid to show the rocket's orientation at various points along the trajectory. Figure \ref{fig: p1 mt states} shows the time histories of the states obtained using the indirect method and DIDO. Figures \ref{fig: p1 mt thrust} shows the thrust magnitude control profile along with its switching function. This time-optimal solution starts at \(T^*=T_\text{min}\) and exhibits one extra thrust magnitude switch than the fuel-optimal solution. Figure \ref{fig: mt gimbal angle} shows that the gimbal angle constraint is saturated for two separate intervals, compared to only one in the fuel-optimal solution. Figure \ref{fig: mt Lagrange gimbal constraint} shows the profiles of \(\tilde{\mu}_4^*\) and \(S_\delta\) accompanying this constraint activation. There seems to be some differences in \(\mu_4\) from each solution, but it evidently does not affect the optimal steering control. Aside from this, all variables of the solution again agree well with those from the DIDO solution, barring discretization errors in the DIDO solution. Inspecting Figures \ref{fig: p1 mt angular velocity magnitude}, \ref{fig: p1 mt glideslope angle}, and \ref{fig: p1 mt tilt angle}, it can be seen the glideslope angle constraint was the only SOPIC that didn't become active. The initial orientation was found to be \(\bm{q}(0)=\left[0.694107905794278, -0.705633365871140, -0.0447276149121824, 0.135259781716865 \right]^\top\). Figures \ref{fig: p1 mt angular velocity magnitude constraint} and \ref{fig: p1 mt tilt angle constraint} show the time histories of the secant penalty functions \(\tilde{S}_{1}\) and \(\tilde{S}_{3}\), respectively, as a function of the tilt angle and angular velocity magnitude, respectively. Figures \ref{fig: p1 mt costates} and \ref{fig: p1 mt costate derivatives} show the time histories of the costate and costate time-derivatives for the time-optimal solution. The same discrepancy (constant offset) in \(\lambda_{\omega_z}\) can be seen in Figure \ref{fig: p1 mt costates}, which is similar to the one reported in Figure \ref{fig: p1 costates}. Figure \ref{fig: p1 mt mass} shows the time histories of the mass, mass flow, mass costate, and mass costate time-derivative. Figure \ref{fig: mt hamiltonian} shows the Hamiltonian time history that satisfied optimality condition ($H(t) = -1$).

\subsection{Analysis of the Lagrange Multipliers associated with State-Only Path Inequality Constraints }

Following Conjecture \ref{conjecture 1}, we calculate the Lagrange multipliers associated with each of the active SOPICs in the fuel- and time-optimal problems and compare the results against DIDO's solutions. We also consider Breakwell problem that has an analytic solution and compare analytic, indirect, and DIDO's solutions. 

\subsubsection{6DOF PDG State-Only Path Constraints Lagrange Multipliers}

\begin{figure}[t!]
    \centering
    \begin{subfigure}[b]{0.49\textwidth}
        \centering
        \includegraphics{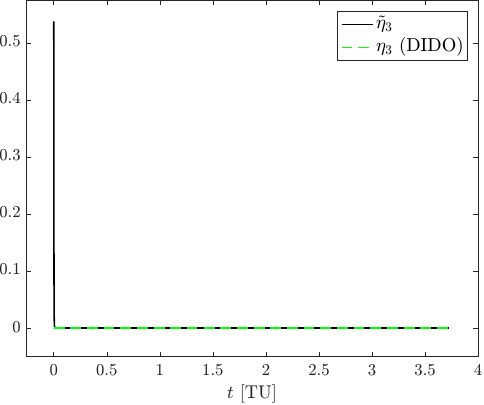}
        \caption{Fuel-Optimal}
        \label{fig: p1 tilt angle constraint multiplier}
    \end{subfigure}%
    ~ 
    \begin{subfigure}[b]{0.49\textwidth}
        \centering
        \includegraphics{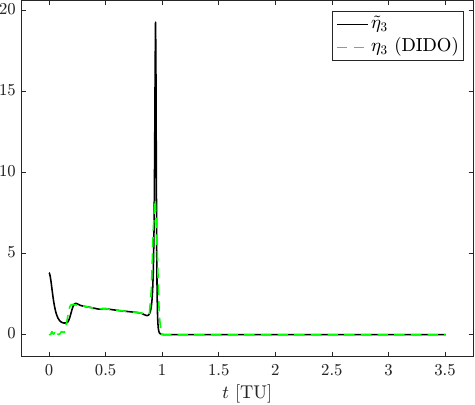}
        \caption{Time-Optimal}
        \label{fig: p1 mt tilt angle constraint multiplier}
    \end{subfigure}
    \caption{Tilt angle constraint multiplier vs. time.}
    \label{fig:comptiltLagrange}
\end{figure}

Figure \ref{fig:comptiltLagrange} shows the Lagrange multiplier associated with the tilt angle constraint for the fuel- and time-optimal solutions. For the fuel-optimal solution, \(\tilde{\eta}_3\), shows an impulse-like spike at the very beginning, which is consistent with the fact that the tilt-angle constraint is a second-order SOPIC (i.e., the control appears after two times differentiation with respect to time). However, DIDO does not show a similar spike. We assume this is because the constraint is only active for an instant at the very beginning of the solution, making it hard for DIDO to capture it (due to the limited number of grid points). On the indirect solution, if we reduced \(\rho_\theta\) even further, then the spike would shrink further into an impulse. For the time-optimal solution, the tilt angle constraint is active for a finite interval. Again, at the beginning of the solution, the indirect solution shows what would likely become an impulse as \(\rho_\theta\) is reduced even further, while DIDO's solution shows no such behavior. While the constraint is active, both solutions agree well, and at the end of the solution, both appear to approximate another impulse where the solution is leaving the constraint arc.

\begin{figure}[t!]
    \centering
    \begin{subfigure}[b]{0.49\textwidth}
        \centering
        \includegraphics{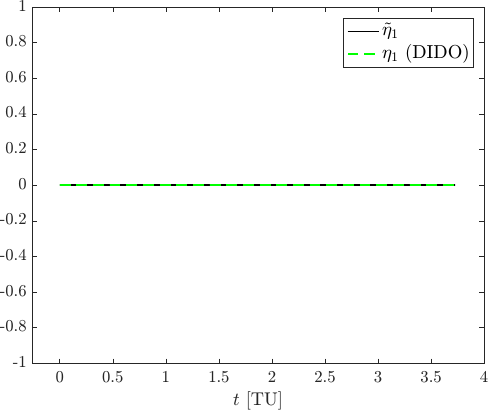}
        \caption{Fuel-Optimal}
        \label{fig: p1 angular velocity magnitude constraint multiplier}
    \end{subfigure}%
    ~ 
    \begin{subfigure}[b]{0.49\textwidth}
        \centering
        \includegraphics{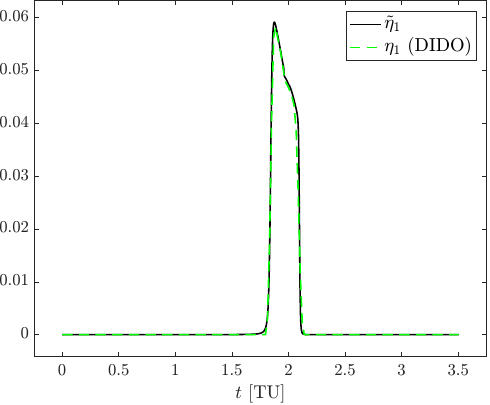}
        \caption{Time-Optimal}
        \label{fig: p1 mt angular velocity magnitude constraint multiplier}
    \end{subfigure}
    \caption{Angular velocity magnitude inequality constraint Lagrange multiplier.}
    \label{fig: angular velocity magnitude constraint multiplier}
\end{figure}

Figure \ref{fig: angular velocity magnitude constraint multiplier} shows the Lagrange multiplier associated with the angular velocity magnitude inequality constraint for the fuel- and time-optimal solutions. The constraint does not become active for the fuel-optimal problem, so \(\tilde{\eta}_1\) and \(\eta_1\) are 0 during the entire maneuver time. For the time-optimal solution, the constraint becomes active for an interval and the behavior of \(\tilde{\eta}_1\) and \(\eta_1\) agree well with each other.

\subsubsection{Breakwell Problem}
As a separate example, we solve the Breakwell problem (Section 3.11 in Ref. \cite{bryson_applied_1975}) with a SOPIC. The Breakwell problem has an analytic solution based on the indirect adjoining approach. Here, we show a comparison between the analytic, indirect, and DIDO's solutions.  This is a minimum-energy problem with a second-order SOPIC, with `$a$' denoting the scalar control input, which is stated as,
\begin{align*}
        \min_{a} \quad & J = \frac{1}{2} \int_{0}^{1} a^2~dt,~ 
        \text{s.t.:}~\dot{x} = v, ~\dot{v} = a,~ x(0)=x(1)=0,~ v(0)=-v(1)=1,~ x(t) \leq l=1/8.
\end{align*}

The problem is solved with the state-path constraint, $S = x(t)-1/8\leq0$, enforced via a secant barrier function augmented to the Hamiltonian and also with DIDO (student version 7.5.6 with 250 nodes). The indirect problem formulation details are omitted for the sake of brevity and only the solutions are presented and compared. 
\begin{figure}[ht!]
    \centering
    \resizebox{\linewidth}{!}{%
        \includegraphics{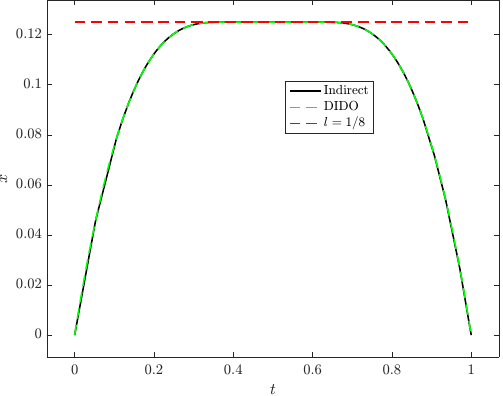}
        \includegraphics{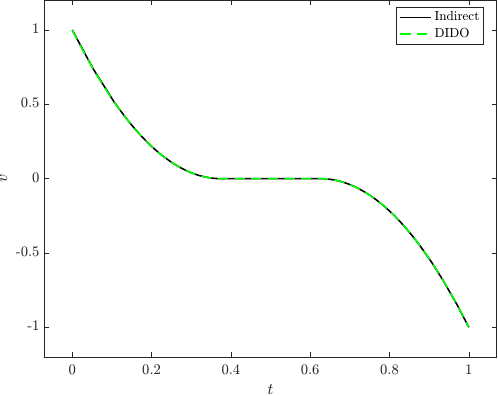}
    }
    \resizebox{\linewidth}{!}{%
        \includegraphics{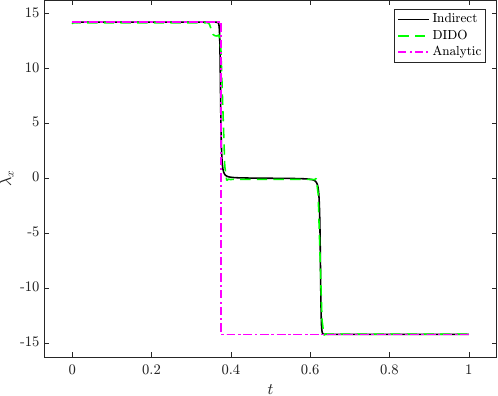}
        \includegraphics{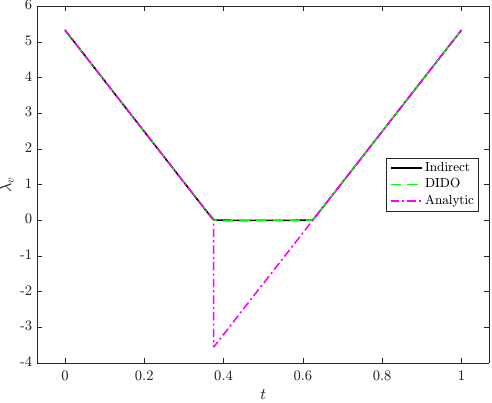}
    }
    \caption{Breakwell problem: states and costates vs. time.}
    \label{fig:breakwellstatescostates}
\end{figure}

The Breakwell problem is a significantly simpler OCP compared to the 6DOF PDG problems. We solved the resulting TPBVP using MATLAB's \verb|bvp5c| with absolute and relative tolerances set to \(1.0 \times 10^{-10}\). Continuation was performed on the penalty function weighting parameter (\(\rho\)) until it was reduced to \(1.0 \times 10^{-10}\). We emphasize that the solutions are obtained independent of each other (i.e., DIDO's solution is not used to initialize the indirect method).  
\begin{figure} [!htb]
    \centering
    \resizebox{0.9\linewidth}{!}{%
        \includegraphics{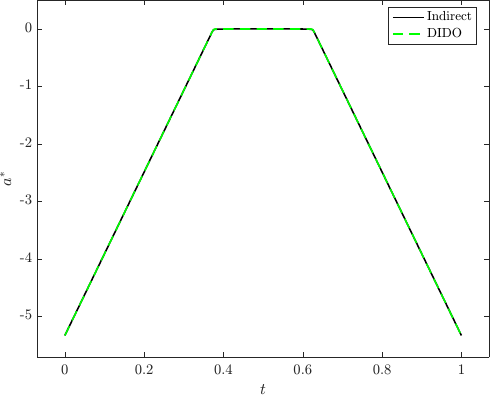}
        \includegraphics{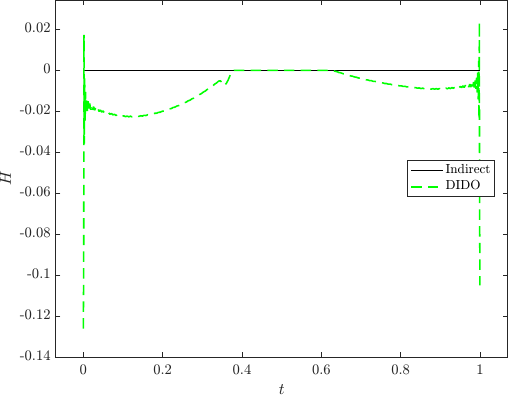}
    }
    \resizebox{0.6\linewidth}{!}{%
        \includegraphics{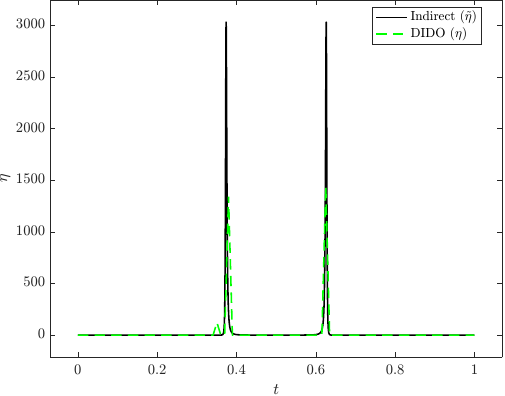}
        
    }
    \caption{Breakwell problem: control, Hamiltonian and constraint Lagrange multiplier vs. time.}
    \label{fig:breakwellcontrolandHamil}
\end{figure}
Figure \ref{fig:breakwellstatescostates} shows the time histories of the states and costates. Figure \ref{fig:breakwellcontrolandHamil} shows the time histories of the control, Hamiltonian, and the Lagrange multiplier associated with the state-path constraint. Upon applying the optimal control theory, we have, $a^* = -\lambda_v$. Both indirect and DIDO's solutions show impulse-type behavior in \(\eta\) at entrance and exit of the constrained arc. This result is consistent with Ref.~\cite{ross_primer_2015} stating that second-order state-path constraints will generally have impulses in the constraint multipliers and jump discontinuities in the costates. The impulses are difficult to numerically capture with DIDO due to only having a finite number of grid points, and difficult to capture with the indirect method due to the approximation of \(\eta\) being limited to the value of \(\rho\). However, the profiles of \(\tilde{\eta}\) and \(\eta\) from DIDO agree well with each other, thus favoring Conjecture \ref{conjecture 1}. The analytic solution gives the same state, control, and Hamiltonian profiles. We only plot the costate profiles in Fig. \ref{fig:breakwellstatescostates} to show the difference in them. Note that, under the implementation approach for the indirect adjoining theory in Ref. \cite{bryson_applied_1975}, the tangency conditions are enforced at the entrance to the constraint arc, resulting in the jump discontinuities in the costates only at entry to the constrained arc.

\section{Conclusion} \label{sec: conclusion}

This work presents an indirect method for solving constrained fuel- and time-optimal six-degree-of-freedom (6DOF) powered descent guidance problems. Inequality constraints are considered on the thrust magnitude, gimbal angle, rocket tilt angle, glideslope angle, and angular velocity magnitude. Derivations of the closed-form control expressions that respect the thrust steering gimbal constraint were presented. State-only path inequality constraints were enforced numerically using interior penalty functions. To overcome the difficulties in solving the multipoint boundary value problems, numerical continuation is used to solve a smooth family of neighboring two-point boundary-value problems. 

Results show that high-accuracy solutions can be obtained that satisfy the necessary conditions of optimality. For the considered problem parameters, boundary conditions, and constraints, the fuel-optimal solution consisted of  2 thrust magnitude switches, whereas the time-optimal solution exhibited 3 thrust magnitude control switches. These problems were also solved with DIDO, which showed similar results in terms of the control, states, and costates profiles. The Lagrange multipliers associated with the direct adjoining of the state-only path constraints to the Hamiltonian were also recovered from the indirect method and DIDO and showed similar trends. Our analysis showed that when no state-only path inequality constraints (SOPICs) were enforced through the penalty function approach, the problem can be easily solved with an indirect single-shooting method. However, the indirect multiple-shooting method was necessary for obtaining solutions under the SOPICs. In fact, on the order of hundreds of shooting segments were required when the SOPICs were active over finite time intervals. This is due to the sensitivity of the problem with respect to the singularity of the secant penalty functions in the neighborhood of the constraint boundaries. The results and comparisons suggest that the proposed indirect method advances the state of the art by applying indirect methods to solve challenging full 6DOF fuel- and time-optimal landing of a rocket in atmosphere subject to state-path nonlinear inequality constraints without requiring a priori knowledge on the structure of the optimal solution. 

\bibliography{refs}

\end{document}